\documentclass[toc,page,title]{article}

\usepackage{PRIMEarxiv}

\usepackage[utf8]{inputenc} 
\usepackage[T1]{fontenc}    
\usepackage{hyperref}       
\usepackage{url}            
\usepackage{booktabs}       
\usepackage{amsfonts}       
\usepackage{nicefrac}       
\usepackage{microtype}      
\usepackage{lipsum}
\usepackage{fancyhdr}       
\usepackage{graphicx}       
\graphicspath{{media/}}     

\usepackage{graphicx}
\usepackage{subcaption}
\usepackage{mathptmx}      
\usepackage[english]{babel}
\usepackage{appendix}
\usepackage{enumitem}
\usepackage{url}            
\usepackage{booktabs}       
\usepackage{amsfonts}       
\usepackage{nicefrac}       
\usepackage{microtype}      
\usepackage{lipsum}
\usepackage{amsthm}
\usepackage{amssymb}
\usepackage{amsmath}
\usepackage{mathptmx}
\usepackage{hyperref}       
\usepackage{xcolor}
\usepackage{comment}
\usepackage{soul}
\usepackage{comment}

\usepackage{algorithm}
\usepackage[noend]{algpseudocode}
\makeatletter
\def\BState{\State\hskip-\ALG@thistlm}
\makeatother

\usepackage[english]{babel}

\DeclareMathOperator*{\argmin}{arg\,min}

\newcommand{\1}{\mathbf{1}}

\newcommand{\ie}{\emph{i.e.}}

 \newtheorem{theorem}{Theorem}
 \newtheorem{lemma}{Lemma}
 \newtheorem{corollary}{Corollary}

 \newtheorem{definition}{Definition}

\title{Dynamic Information Manipulation Game
\thanks{This version: July 2025.} 
}

\author{
  Shutian Liu,\textsuperscript{a} Quanyan Zhu,\textsuperscript{b} \\
 \textsuperscript{a}Department of Systems Engineering, City University of Hong Kong, Hong Kong, China;\\
\textsuperscript{b}Tandon School of Engineering, New York University, Brooklyn, NY 11201, USA  \\
  \texttt{Contact: shutian.liu@cityu.edu.hk(SL); qz494@nyu.edu(QZ)} \\
}

\begin{document}
\maketitle

\begin{abstract}
We propose a dynamic information manipulation game (DIMG) to investigate the incentives of an information manipulator (IM) to influence the transition rules of a partially observable Markov decision process (POMDP). DIMG is a hierarchical game where the upper-level IM stealthily designs the POMDP's joint state distributions to influence the lower-level controller's actions. DIMG's fundamental feature is characterized by a stagewise constraint that ensures the consistency between the unobservable marginals of the manipulated and the original kernels. In an equilibrium of information distortion, the IM minimizes cumulative cost that depends on the controller's informationally manipulated actions generated by the optimal policy to the POMDP. 
We discuss ex ante and interim manipulation schemes and show their connections.
The effect of manipulation on the performance of control policies is analyzed through its influence on belief distortion.
\end{abstract}


\section{Introduction}
\label{sec:intro}
Decisions are often made without complete knowledge of the underlying environment.
In addition to this limitation, decision-makers (DMs) also have to take into account the credibility of the available information, since recent advances in social media, information technologies, and artificial intelligence have made the generation and distribution of misinformation convenient and inexpensive.
With the actions and operations performed in modern complex systems becoming increasingly dynamic and strategic, information about either the states of the environment or the preferences of others that are hardly credible can significantly influence the choices of DMs in a stealthy and silent way.
Actions based on manipulated information lead, for example, to security vulnerabilities in networked cyber systems, cascading operational failures in industrial supply chains, and fatal accidents in autonomous vehicle systems.
Therefore, it is crucial to understand the misinformation generation procedure and the potential impact of it on decision making and system performance.

There is a vast literature on sequential decision-making that involves the modeling of imperfect information from the DM's perspective.
One popular approach suggests introducing uncertainty to the DM's model parameters.
Robust controllers are obtained under the worst-case scenario with respect to certain criteria in the classic monographs (see, e.g., \cite{bacsar2008h,khalil1996robust}).
A generalization of the robustness approach focuses on the setting where the distributional information of the model parameters is unknown, which is referred to as ambiguity in the economics literature (e.g., \cite{epstein2003recursive,maccheroni2006ambiguity}).
For sequential decision-making problems, this includes the robust dynamic programming, e.g., \cite{iyengar2005robust,nilim2005robust,wiesemann2013robust}, where the transition kernel of a Markov decision process (MDP) is subject to ambiguity, and the works that pertain to the topics of distributionally robust MDP (see \cite{xu2010distributionally,nakao2021distributionally}), distributionally robust control (e.g., \cite{van2015distributionally,yang2020wasserstein}), 
recursive robust estimation and control (e.g., \cite{hansen2005robust,hansen2007recursive}), and the ambiguous partially observable Markov decision process (POMDP) in \cite{saghafian2018ambiguous}.
A closely related approach that emphasizes the preference aspect of DMs is referred to as dynamic risk-sensitive stochastic optimization (e.g., \cite{ruszczynski2010risk,pflug2014multistage,bauerle2022distributionally}).
Risk sensitivity has been shown to correspond to ambiguity in the transition laws of the state process.
Many of the above works have also emphasized the nature that robust (and distributionally robust) decision-making frameworks can be equivalently formulated as zero-sum games between the DM and an adversary who controls the unknown parameter; see, for example,  \cite{bacsar2008h,ruszczynski2010risk,james1994risk,bauerle2022distributionally,hansen2005robust,saghafian2018ambiguous}.

Uncertainties and ambiguities are commonly described as predefined parameter sets either constructed based on experience or estimated from data generated by nominal models. 
However, the equivalent formulation of robustness in the form of an adversary in zero-sum games suggests that, apart from using fixed parameter sets to describe the unknowns, one can also regard the sources of the unknowns as other decision-making parties and study their interactions with the DM.  
This type of relations is frequently observed in practice in scenarios such as data poisoning attacks in recommendation systems (e.g., \cite{zhang2021data}), misguidance in traffic equilibria (e.g., \cite{pan2022poisoned}), and intentional information dissemination in social media (e.g., \cite{wu2024attacking}).
Furthermore, misinformation shapes individual behaviors or even populational behavioral patterns in a stealthy manner.
As documented in \cite{sharma2019combating} and \cite{abdali2022multi}, the identification of misinformation is itself a very challenging task.
Meanwhile, timely decisions are often made before one can gather an extensive amount of information.
Local and partial information are naturally biased, and they can even be intentionally filtered.
These facts motivate us to take a step backward to reconsider the generation and impact of misinformation on decision-making and re-examine the level of informational advantage of a malicious information producer.

In this paper, we take an alternative approach and investigate the creation and impact of misinformation from the perspective of an information manipulator (IM).
We propose a dynamic information manipulation game (DIMG) based on the generic risk-sensitive POMDP framework in \cite{baauerle2017partially} and introduce the IM as a designer who presents the DM with a simulacrum of the stochastic kernel of the POMDP.
In particular, the IM is in place of the stochastic kernel and designs the joint distribution of the observable and the unobservable states at each stage.
In the meantime, the DM receives samples of the observable states from the IM thinking that they are generated from the true kernel and optimizes her objective function as if the IM is absent.
This setting fits the property that one sample of the observable state per stage is sufficient for the DM to decide her control action, yet the sample is insufficient in justifying its own credibility.
The IM's designs are chosen to minimize his own cumulative cost that depends on the actions of the DM.
Referred to as an information manipulator, the IM's power is limited when it comes to the unobservable part of the state. 
This restriction is setup by introducing a stage-wise consistency constraint to the IM's manipulation problem, which requires the IM to adhere to the marginal distribution of the unobservable state specified by the true kernel given history.
In other words, the IM can only manipulate the joint distribution up to a given marginal distribution of the unobservable state at each stage.

\textcolor{black}{
We introduce two manipulation schemes available to the IM that correspond to different information patterns.
The first scheme is called ex ante manipulation scheme, where the IM chooses the joint probability distribution of the observable and unobservable states before the realization of the true kernel at each stage.
The second one is referred to as interim manipulation scheme, where the IM acts after he observes the realization of the kernel at each stage.
In this scheme, determining the joint distribution reduces to identifying the distribution of the observable state conditioned on the realization of the unobservable state.
}
DIMG is a hierarchical game with the IM in the upper level and the DM in the lower level.
However, DIMG differs from standard dynamic Stackelberg games since the lower-level player only observes a random sample rather than the whole distribution of information designed by the upper-level player.
\textcolor{black}{
In both schemes, the DM's available information is contained in set of available information of the IM at each stage.
Hence, the IM is able to reconstruct the DM's optimal control policy.
We introduce equilibrium of information distortion in the both schemes as the solution concept to DIMG.
At an equilibrium of information distortion, the IM designs the state distributions feasible to the chosen scheme so that the DM's control actions generated by the optimal policy to the POMDP given observations following the IM's designs optimize the IM's objective function.
}
As we will show later that, in DIMG, the DM's actions depend on the IM's action only through the observation samples and the DM's policies remain unchanged.

The analysis of the DIMG is divided into the following three steps.
We first follow \cite{baauerle2017partially} and consider the variant where unnormalized conditional densities of the unobservable state and accumulated cost given observable history are adopted as information states to reformulate the risk-sensitive POMDP with a generic utility function into a completely observable MDP.
While the approach introduced in \cite{baauerle2017partially} has advanced the techniques of \cite{james1994risk} and \cite{di1999risk} where partially observed stochastic control problems with exponential utility are solved by including the stage costs into the information states, we have addressed the gap where \cite{baauerle2017partially} only considers normalized conditional densities.
The unnormalized conditional densities help simplify the dynamic programming recursions and the analysis of the impacts of information manipulation in the later parts.
After introducing the dynamic programming equations of the DM, the second part of our analysis incorporates the IM.
We zoom in on the optimization problem of the IM at each stage.
The existence of optimal designs of the IM and the connection between the solutions under the two schemes are shown using a technique built on the ``joint$+$marginal" approach introduced in \cite{lasserre2010joint+}.
In fact, the stage-wise consistency constraints, which specify marginal distributions, propagate forward in time and the joint distribution design is solved backward in time based on the dynamic programming equations of the DM.
\textcolor{black}{
Finally, in the third part, we investigate the impact of information manipulation on the performance of the DM's control policy.
This is achieved by leveraging the represention of the DM's objective function under the information state MDP and analyzing the influence of informationally distorted beliefs.
We derive a tight upper-bound on the difference between the performances of the DM's policy when the IM is present and absent.
}

The remainder of the paper is organized as follows.
In Section \ref{sec:related}, we briefly summarize additional related studies.
We introduce the preliminary POMDP model and the DIMG in Section \ref{sec:model}.
In Section \ref{sec:solution}, we first solve the DM's problem using unnormalized conditional densities.
Then, we derive the solution to the DIMG.
Section \ref{sec:impacts} investigates the impact of information manipulation on the performance of the DM's problem.
Section \ref{sec:conclusion} concludes the paper and suggests potential applications.

\section{Related Works}
\label{sec:related}

Information manipulation belongs to the broad field of adversarial attacks in dynamic decision making, which includes, but not limited to, perturbation analysis in stochastic control problems (e.g., \cite{kushner2012weak,borkar2010singular}), optimal design of adversarial attacks in MDPs (e.g., \cite{russo2019optimal,mcmahan2024optimal}), and reinforcement learning with poisoning attacks (e.g., \cite{ma2019policy,liu2021provably,do2023adversarial}).
In particular, \cite{mo2010false} and \cite{liu2011false}, among the others, consider data injection attacks to model the malicious incentives for controlling system measurements.
Information manipulation differentiates from data injection attacks since the IM simulates the stochastic kernel of the POMDP (hence, the manipulated kernel is referred to as a simulacrum of the true kernel) instead of merely designs the observations.
The stealthiness of the IM is ensured by the construction of the DIMG, while data injection attacks need to be designed so that they do not trigger the alarms from the failure detection devices (see \cite{mo2010false}).
We use the IM in the DIMG to model intentional third parties supported by advanced information technology (IT), who aims for stealthy information manipulation in the long-term.

Our work is also closely related to the recent work of \cite{zheng2022privacy}, where the authors consider belief manipulation in POMDPs.
Zero-sum stochastic games are commonly used to model the relations between adversaries and DMs, see for example \cite{bauerle2017zero} and \cite{nguyen2009stochastic}.
Among them, frameworks built on stochastic games with one-sided partial observations, such as \cite{horner2010markov,chatterjee2014partial,zheng2022stackelberg}, and \cite{horak2023solving}, are closely related to our framework.
In particular, \cite{zheng2022stackelberg} has studied Stackelberg equilibria with an informationally inferior follower.

The design of information from a rational designer's perspective has been extensively studied in the economic literature. 
Classic approaches include, for example, the signaling games of \cite{spence1978job}, the cheap talks of \cite{crawford1982strategic}, and the verifiable messages of \cite{grossman1981informational}.
Recently, there has been an increasing interest in Bayesian persuasion (see \cite{kamenica2011bayesian}), which is also referred to as information design (e.g., \cite{bergemann2019information,kamenica2019bayesian}). 
The above approaches investigate the design of information in a setting where the information receiver is aware of the rule of information generation and responds to observations rationally.
In DIMG, we further strengthen the informational advantage of the designer by not forcing the revealing of the designed information.
Meanwhile, the designer's power is also restricted by the stage-wise consistency constraint.

The information state recursion, which is also referred to as belief updates or Bayesian filter, has been widely used to solve sequential decision-making problems when the DM only has partial observations of the system, see \cite{bertsekas1996stochastic,kumar2015stochastic,bensoussan1992stochastic,bauerle2011markov,elliott2008hidden,feinberg2016partially}, and the references therein.
While essentially the information states are conditional probability distributions of the unobservable states given the current observable history, many variations on this basic setting have been investigated.
One variation adopts unnormalized conditional densities as information states, see, for example, \cite{kumar2015stochastic,elliott2017discrete}, and \cite{bensoussan1992stochastic}.
Another extension integrates stage costs into information states to address risk-sensitive objectives, including \cite{james1994risk} and \cite{di1999risk}, where the setting of exponential cost and a change of measure technique are adopted and the recent work \cite{baauerle2017partially} where the authors' derivation is based on a general POMDP setting without requiring the change of measure.
There are also works, such as \cite{subramanian2022approximate,dave2023approximate}, and \cite{yang2022discrete}, which have adopted approximate information states when the system is not exactly known.

Finally, the analysis of the impact of information manipulation is related to filter stability results in stochastic control problems or POMDPs in the literature, see \cite{le2004stability,kara2020robustness,mcdonald2022robustness,kara2022robustness,golowich2023planning}, and the reference therein.
Similar results are also discussed in \cite{saghafian2018ambiguous} under the ambiguous POMDP setting.
The recent work \cite{zheng2022provable} has emphasized the impacts of certain special observation distributions on the performance of a POMDP, which is closely related to our focus on the impacts of the IM.

\section{Dynamic Information Manipulation Game}
\label{sec:model}
In this section, we introduce DIMG.
We first present the DM's POMDP at the lower level and then motivate the IM's design problem at the upper level.

\subsection{Lower Level Risk-Sensitive POMDP}
\label{sec:model:pomdp}
\textcolor{black}{Consider a process $(X_n, Y_n)_{n\in \mathbb{N}_0}$ on the Borel state space $\mathcal{X}\times \mathcal{Y}$ with finite time horizon $N$. 
}
The process $X_n$ is observable to the DM, while the process $Y_n$ is unobservable.
Let the Borel space $\mathcal{A}$ denote the set of all available actions of the DM. 
The subset $\mathcal{D}\subset \mathcal{X}\times \mathcal{A}$ contains graph of a measurable mapping from $\mathcal{X}$ to $\mathcal{A}$.
We use $\mathcal{D}(x):=\{a\in \mathcal{A}:(x,a)\in \mathcal{D}\}$ to denote the feasible action given the observable state $x\in\mathcal{X}$.
We assume that $\mathcal{D}(x)$ is convex and compact for all $x\in\mathcal{X}$.
The transition of the states is described by the stochastic kernel $Q: \mathcal{D}\times \mathcal{Y} \rightarrow \mathcal{P} (\mathcal{X}\times \mathcal{Y})$, where $\mathcal{P} (\mathcal{X}\times \mathcal{Y})$ denotes the set of probability measures on $\mathcal{X}\times \mathcal{Y}$.
This means that, given the current pair $ (x,y)$ and an action $a\in \mathcal{D}(x)$, the probability that the next state pair is in a Borel set $B\in \mathcal{B}(\mathcal{X}\times \mathcal{Y})$ is $Q(B|x,y,a)$.
Assume that the stochastic kernel $Q$ has a measurable density function $q$ with respect to the reference probability measures $\lambda$ on $\mathcal{P}(\mathcal{X})$ and $\nu$ on $\mathcal{P}(\mathcal{Y})$, \ie, 
\begin{equation*}
    Q(B|x,y,a)=\int_{B}q(x',y'|x,y,a)\lambda(dx')\nu(dy'), \forall B\in \mathcal{B}(\mathcal{X}\times \mathcal{Y}).
\end{equation*}
We assume that $q$ is continuous and bounded with respect to all of its arguments.
\textcolor{black}{
The $X$-marginal and the $Y$-marginal of the transition kernel density $q$ are denoted by $q^X(x'|x,y,a):=\int_{\mathcal{Y}}q(x',y'|x,y,a)\nu(dy')$ and $q^Y(y'|x,y,a):=\int_{\mathcal{X}}q(x',y'|x,y,a)\lambda(dx')$, respectively.
The $X$-marginal and the $Y$-marginal of transition kernel $Q$ are denoted by $Q^X(B|x,y,a):=\int_{B}q^X(x'|x,y,a)\lambda(dx')$ for $B\in\mathcal{B}(\mathcal{X})$, and $Q^Y(B|x,y,a):=\int_{B}q^Y(y'|x,y,a)\nu(dy')$ for $B\in\mathcal{B}(\mathcal{Y})$, respectively.
}
We also assume that the initial distribution $Q_0^Y$ of $Y_0$ is given and $Q_0^Y$ has measurable density $q_0^Y$ with respect to $\nu$.

The observable histories of the DM are defined as:
\begin{equation*}
    \begin{aligned}
        & H_0:= \mathcal{X}, \\
        & H_n:= H_{n-1} \times \mathcal{A} \times \mathcal{X}.
    \end{aligned}
\end{equation*}
An element of $H_n$ is denoted $h_n:=(x_0,a_0,x_1,\cdots,a_{n-1},x_n)$.
The control policy of the DM is denoted as $\pi:=(g_0,g_1,\cdots)\in \Pi$, where $g_n:H_n\rightarrow \mathcal{D}(x_n)$ is the strategy at stage $n$.
Under a policy $\pi\in\Pi$, the sequence of actions chosen is defined as 
\begin{equation*}
    \begin{aligned}
        &A_0:=g_0(X_0),\\
        &A_n:=g_n(X_0,A_0,X_1,\cdots,A_{n-1},X_n).
    \end{aligned}
\end{equation*}

Given a policy $\pi$ and an initial state pair $(X_0,Y_0)=(x,y)$, the transition kernel $Q$ determines a probability measure $\mathbb{P}^\pi_{xy}$ on the product $\sigma$-algebra on $(\mathcal{X}\times\mathcal{Y})^{N+1}$ according to the Ionescu-Tulcea theorem.
In the scenario where only the distribution $Q_0^Y$ of $Y_0$ is known, the measure is denoted as $\mathbb{P}^\pi_x:=\int_{\mathcal{Y}}\mathbb{P}^\pi_{xy}Q_0^Y(dy)$.

The stage cost received by the DM is represented by a continuous and bounded function $c:\mathcal{X}\times \mathcal{Y} \times \mathcal{A}\rightarrow \mathbb{R}_+$.
The DM's utility function $U:\mathbb{R}_+\rightarrow \mathbb{R}$ is assumed continuous and strictly increasing.
\textcolor{black}{
A utility function represents the DM's risk attitude towards uncertain outcomes in risky environments.
When considering a minimization problem, the DM is risk-averse if $U$ is convex and risk-seeking if $U$ is concave.
}
Given a policy $\pi$ and an initial observation $X_0=x$, the DM's objective function can be formulated as:
\begin{equation}
   J_{N\pi}(x):=\int_{\mathcal{Y}}\mathbb{E}^{\pi}_{xy}\left[ U\left( \sum_{k=0}^{N-1}\beta^k c(X_k,Y_k,A_k) \right) \right] Q_0^Y(dy),
   \label{eq:obj func}
\end{equation}
where $0<\beta<1$ is a discount factor.
The DM's optimization problem is $J_N(x):=  \inf_{\pi} J_{N\pi}(x)$.
\textcolor{black}{
Note that considering $J_N(x)$ is equivalent to considering minimizing the certainty equivalent of the accumulated costs, which is defined as $U^{-1}(J_N(x))$ in our setting.
The reason lies in that the utility function $U$ is strictly increasing.
We refer the reader to \cite{bauerle2011markov} for more discussion of utility theory and its use in MDPs.
}

\subsection{Upper Level Information Manipulation Problem}
\label{sec:model:information manipulation}
The approach to optimizing (\ref{eq:obj func}) from the DM's perspective will be discussed in detail in Section \ref{sec:solution}.
Nevertheless, we briefly recall the procedure of solving a generic POMDP here to motivate the DIMG.

A common technique to solve a POMDP is to transform it into a completely observable MDP with the help of information states.
The information states, which are recursively constructed by computing the conditional probability distribution of the unobservable states given the current observation using the Bayes operator, serve as a ``sufficient statistic" for representing the DM's objective function.
In fact, one sample of the observable state is enough for the DM to update her information state at each stage.
On the one hand, this sufficiency aids the DM in solving her control problem with the minimal amount of information.
On the other hand, the information sufficient for controlling the process is insufficient to justify the credibility of itself.
In other words, it is challenging for the DM to examine whether the sampled observation is drawn from the true transition kernel or not.
This issue creates space for a third party to stand in the middle of the DM and the kernel, and conduct information manipulation.

Consider the variation of the setting introduced in Section \ref{sec:model:pomdp} where each stage is divided into two levels.
\textcolor{black}{
At the higher level, a stealthy IM designs the joint probability measure $P_n(x,y)\in\mathcal{P}(\mathcal{X}\times\mathcal{Y})$ of the state $(X_n, Y_n)$ of stage $n=0,1,\cdots,N-1$.
}
We assume that $P_n$ has a measurable density $p_n$ with respect to the probability measures $\lambda$ and $\nu$, \ie,
\begin{equation*}
    P_n(B)=\int_{B}p_n(x,y)\lambda(dx)\nu(dy), \  \ \forall B\in\mathcal{B}(\mathcal{X}\times\mathcal{Y}).
\end{equation*}
\textcolor{black}{
Let $p_n^X(x):=\int_{\mathcal{Y}}p_n(x,y)\nu(dy)$ and $p_n^Y(y):=\int_{\mathcal{X}}p_n(x,y)\lambda(dx)$ denote the $X$-marginal  and the $Y$-marginal of the density $p_n$, respectively.
Let $P_n^X(B):=\int_{B}p_n^X(x)\lambda(dx)$ for $B\in\mathcal{B}(\mathcal{X})$ and $P_n^Y(B):=\int_{B}p_n^Y(y)\nu(dy)$ for $B\in\mathcal{B}(\mathcal{Y})$ denote the $X$-marginal and the $Y$-marginal of the $P_n$, respectively.
}
At the lower level, the state $(x_n, y_n)$ realizes and follows the distribution $P_n$.
The DM then chooses her control after receiving the observable part $x_n$.
We refer to the IM as a designer who provides a simulacrum of the stochastic kernel, as he assumes the role of the stochastic kernel, with the DM making decisions based on the information he supplies.
The stealthiness of the IM is justified by the assumption that the DM is unable to figure out whether the sampled observations come from the simulacrum or the kernel itself.
In other words, while the DM solves her control problem thinking that the observations are sampled from the true distribution of the transition kernel as if the IM is absent, the information that she actually receives follows the manipulated distributions chosen by the IM.
Note that this assumption fits the property that one sample of observation at each stage is sufficient for decision-making.
We restrict the IM's power so that he only has control over the observable part of the state process, justifying the name information manipulator.
This setting is fulfilled by involving a stage-wise consistency constraint of the feasible designs which requires the $Y$-marginal of the designed joint distribution to be equal to the $Y$-marginal of the joint distribution obtained by the kernel $Q$ given the previous state and action.
The definition of a stagewise consistent design is presented as follows.
\textcolor{black}{
\begin{definition}
\label{def:stagewise consistency}
(Stagewise Consistency.)
The IM's choice of joint state distribution $P_n$ is called stagewise consistent, if for $n=1,\cdots,N-1$ and for given $x_{n-1},y_{n-1},a_{n-1}$, the following equality holds:
\begin{equation}
    P^Y_n(B)=Q^Y(B|x_{n-1},y_{n-1},a_{n-1}), \  \ \forall B\in\mathcal{B}(\mathcal{Y}),
    \label{eq:consistency 1 to N}
\end{equation}
and for $n=0$, it holds that
\begin{equation}
    P^Y_0(B)=Q_0^Y(B),\  \ \forall B\in\mathcal{B}(\mathcal{Y}).
    \label{eq:consistency 0}
\end{equation}
\end{definition}
Let $\theta:=(P_0,P_2,\cdots,P_{N-1})\in \Theta\subset [\mathcal{P}(\mathcal{X}\times\mathcal{Y})]^N$ denote the concatenation of $N$ designs, where we use $\Theta$ to denote the set of designs that are stagewise consistent.
}
Note that, since the IM cannot directly influence the probability laws of the unobservable part of the state process, it is the joint effort of both the IM and the DM that governs the evolution of the process $(X_n,Y_n)$.

While the DM can only observe a sample of the observable state that follows the distribution $P_n$ at the lower level of each stage, we consider two possible design schemes corresponding to different information patterns at the higher level of each stage.
The first one is called the ex ante manipulation scheme, where the IM determines $P_n$ without observing the realization from the true kernel at stage $n$.
Note that the IM does recall the history up to stage $n-1$, which consists of both the realizations of the observable and the unobservable states, \ie,  $(x_0,y_0,x_1,y_1,\cdots,x_{n-1},y_{n-1})$, the actions of the DM, \ie,  $(a_0,a_1,\cdots,a_{n-1})$, and the designs of the IM, \ie,  $(P_0,P_1,\cdots, P_{n-1})$.
The second one is referred to as the interim manipulation scheme, where the IM first observe the realization of the true states $(x_n, y_n)$ from the kernel at stage $n$ in addition to the observations available under the ex ante scheme, then he chooses $P_n$ in the collapsed version, \ie, the distribution of the observable part given the realization $y_n$. 
It is straightforward that interim information manipulation is a special case of ex ante information manipulation.
Their connections will be elaborated in Section \ref{sec:solution:solve info manipulation}.
For now, we focus on the ex ante scheme unless specified otherwise.

Let $\mathbb{P}^{\pi\theta}$ denote the probability measure on the product $\sigma$-algebra on $(\mathcal{X}\times\mathcal{Y})^{N+1}$ generated by the design $\theta=(P_0,\cdots,P_{N-1})\in\Theta$ given initial distribution of the unobservable state denoted $Q_0^Y$.
\textcolor{black}{
The reason why the measure $\mathbb{P}^{\pi\theta}$ is dependent on both $\theta$ and $\pi$ is because of stagewise consistency.
}
Though we do not make $Q_0^Y$ explicit in $\mathbb{P}^{\pi\theta}$ for notational simplicity, it is clear from (\ref{eq:consistency 0}) that the measure $\mathbb{P}^{\pi\theta}$ depends on the initial distribution through $\theta$.
The IM's optimization problem is $I_{N\pi}:=\inf_{\theta\in\Theta} I_{N\pi\theta}$ driven by the following objective function
\begin{equation}
    I_{N\pi\theta}:=\mathbb{E}^{\pi\theta}\left[ \sum_{k=0}^{N-1}\alpha^k \left( r(X_k,Y_k,A_k)+ \rho_k(P_k) \right) \right] ,
    \label{eq:designer's problem ex ante}
\end{equation}
where $r:\mathcal{X}\times\mathcal{Y}\times\mathcal{A}\rightarrow\mathbb{R}_+$ denotes a continuous and bounded cost function dependent on the states and the DM's action and $\rho_k:\mathcal{P}(\mathcal{X}\times \mathcal{Y})\rightarrow\mathbb{R}_+$ captures the cost of performing information manipulation at stage $k$. 
\textcolor{black}{
A scalar $\gamma\geq0$ can be introduced to balance these two costs in (\ref{eq:designer's problem ex ante}).
For notational convenience, we assume that it has been absorbed into $r$.
}
Meaningful choices for the function $\rho_k$ include well-defined functions that measure the deviation from the designed joint distribution $P_n$ to the joint distribution obtained from the kernel $Q$ given previous states and action.
One such function is the total variation distance between the two distributions defined as $d_{TV}(P,Q):=\frac{1}{2}\int|dP-dQ|$.
In this paper, we let $\rho_n(P_n)=2d_{TV}(P_n,Q(\cdot,\cdot|x_{n-1},y_{n-1},a_{n-1}))$ for notational convenience.

The solution concept of the DIMG is presented in the following.
\begin{definition}
(Ex Ante Equilibrium of Information Distortion.)
In the DIMG with an upper-level IM and a lower-level DM under the ex ante information pattern, the IM's design of joint probability distributions $\theta^*:=(P_0^*,P_2^*,\cdots,P_{N-1}^*)\in\Theta$  and the DM's control policy $\pi^*:=(g_0^*,\cdots,g_{N-1}^*)\in\Pi$ constitute an ex ante equilibrium of information distortion, if $I_{N\pi^*\theta^*}\leq I_{N\pi^*\theta}$ for all stagewise consistent designs $\theta \in \Theta\subset[\mathcal{P}(\mathcal{X}\times\mathcal{Y})]^N$ and $J_{N\pi^*}\leq J_{N\pi}$ for all $\pi\in\Pi$.
The quantity $I_{N\pi^*\theta^*}$ is the IM's ex ante manipulation cost.
\end{definition}

The concept of interim equilibrium of information distortion can be defined in a similar manner if we restrict the designs of the IM to be the probability distribution of the observable state given the realization of the current unobservable state.

\section{Analysis of Dynamic Information Manipulation Game}
\label{sec:solution}
This section focuses on the analysis of DIMG.
We first consider the DM's POMDP ignoring the IM and find its optimal policy.
In particular, we show in Section \ref{sec:solution:info states} that unnormalized conditional density functions of the unobservable states and the accumulated costs serve as the information states and in Section \ref{sec:solution:DP equations} how these information states help construct an equivalent MDP with extended state space.
We derive the dynamic programming equations of this MDP.
Then, in Section \ref{sec:solution:solve info manipulation}, we involve the IM and derive the solution to the DIMG using the DM's dynamic programming equations as building blocks.
The connection between the ex ante and the interim equilibria of information distortion is also discussed.
Two examples are considered in Section \ref{sec:example} to illustrate the procedure of information manipulation, to present tractable reformulations of the IM's problem for computational purposes, and to generate insights about the maneuverability of information.

\subsection{Decision-Maker's Information States}
\label{sec:solution:info states}
Let $\mathcal{M}(\mathcal{Y}\times \mathbb{R}_+)$ denote the set of measures on the $\sigma$-algebra $\mathcal{B}(\mathcal{Y}\times \mathbb{R}_+)$.
We consider $\mu\in\mathcal{M}(\mathcal{Y}\times \mathbb{R}_+)$ as an information state for the partially observed model.
The measure $\mu$ encodes the unnormalized conditional distribution of the hidden state and the accumulated cost.
Let $\mu^Y(dy):=\int_{\mathbb{R}_+}\mu(dy,ds)$ and $\mu^S(ds):=\int_{\mathcal{Y}}\mu(dy,ds)$ denote the $Y$-marginal and the $S$-marginal of an information state $\mu$, respectively.

For $B\in \mathcal{B}(\mathcal{Y}\times \mathbb{R}_+)$, consider the update rule $\Psi:\mathcal{X}\times\mathcal{A}\times\mathcal{X}\times\mathcal{M}(\mathcal{Y}\times \mathbb{R}_+)\times\mathbb{R}_+\rightarrow \mathcal{M}(\mathcal{Y}\times \mathbb{R}_+)$ defined by
\begin{equation}
    \Psi(x,a,x',\mu,z)(B):=\int_{\mathcal{Y}}\int_{\mathbb{R}_+}\left( \int_{B}q(x',y'|x,y,a)\nu(dy')\delta_{s+zc(x,y,a)}(ds') \right) \mu(dy,ds),
    \label{eq:info state operator}
\end{equation}
where $\delta_s(\cdot)$ denotes the Dirac measure at point $s\in\mathbb{R_+}$.
\textcolor{black}{
Define for $B\in\mathcal{B}(\mathcal{Y}\times \mathbb{R}_+)$, and $h_n$ the following sequence of measures
\begin{equation}
    \begin{aligned}
        \mu_0(B|h_0)&:=(Q_0^Y\times\delta_0)(B), \\
        \mu_{n+1}(B|h_n,a,x')&:=\Psi(x_n,a,x',\mu_n(\cdot|h_n),\beta^n)(B), n=0,\cdots,N-1.
        \label{eq:info state recursion}
    \end{aligned}
\end{equation}
Let $D_n:=\int_{\mathcal{Y}}q^X(x_n|x_{n-1}, y_{n-1}, a_{n-1})\mu_{n-1}^Y(dy_{n-1})$ for $n=1,2,\cdots, N$, and $D_0:=1$. 
We observe that $D_n$ is the normalization constant of (\ref{eq:info state operator}) for stage $n$.
Define the following random variables
\begin{equation*}
\begin{aligned}
    & S_0:=0, \\
    & S_n:=\sum_{k=0}^{n-1}\beta^kc(X_k,Y_k,A_k), n=1,\cdots,N.
    \end{aligned}
\end{equation*}
The next result shows the role of the sequence $(\mu_n)$.
}
\textcolor{black}{
\begin{theorem}
 Suppose that the information state process $(\mu_n)$ is defined by (\ref{eq:info state recursion}). 
 Then, it holds for $n=0,\cdots,N$, $\pi\in\Pi$, and $B\in\mathcal{B}(\mathcal{Y}\times\mathbb{R}_+)$ that
 \begin{equation*}
     \prod_{k=0}^nD_k\cdot\mathbb{P}_x^\pi((Y_n,S_n)\in B|X_0,A_0,\cdots,X_n)=\mu_n(B|X_0,A_0,\cdots,X_n),\  \  \mathbb{P}_x^\pi\text{-a.s.}.
 \end{equation*}
\end{theorem}
}
\textcolor{black}{
\noindent{\bf Proof}.
We first show for $n=0,\cdots,N$ that the relation
\begin{equation}
    \prod_{k=0}^n D_k \cdot \mathbb{E}_x^{\pi}[v(X_0,A_o,X_1,\cdots,X_n,Y_n,S_n)]
    = 
    \mathbb{E}_x^{\pi}[v'(X_0,A_0,X_1,\cdots,X_n)]
    \label{eq:condition for info state recursion}
\end{equation}
holds for all bounded and measurable $v:H_n\times \mathcal{Y}\times \mathbb{R}_+\rightarrow \mathbb{R}$ and 
$v'(h_n):=\int_{\mathcal{Y}}\int_{\mathbb{R}_+}v(h_n, y_n, s_n)\mu_n(dy_n, ds_n|h_n)$ using induction.
For $n=0$, both sides of (\ref{eq:condition for info state recursion}) reduce to $\int_{\mathcal{Y}}v(x_0,y_0,0)Q_0^Y(dy)$.
Suppose that (\ref{eq:condition for info state recursion}) holds for $n-1$.
Then, given $h_{n-1}$, the left-hand side of (\ref{eq:condition for info state recursion}) becomes
\begin{equation*}
    \begin{aligned}
        \prod_{k=0}^n D_k\cdot \mathbb{E}_x^{\pi}[v(h_{n-1},A_{n-1}, X_n, Y_n, S_n)]
        =& D_n \cdot \int_{\mathcal{Y}}\int_{\mathbb{R}_+}\mu_{n-1}(dy_{n-1},ds_{n-1}|h_{n-1})\\
        & \  \ \cdot \int_{\mathcal{Y}}\int_{\mathbb{R}_+}\nu(dy_n)\lambda(dx_n)q(x_n,y_n|x_{n-1},y_{n-1},g_{n-1})\\
        & \  \ \cdot
        \int_{\mathbb{R}_+}\delta_{s_{n-1}+\beta^{n-1}c(x_{n-1},y_{n-1},g_{n-1})}(ds_n)v(h_{h-1},g_{n-1},x_n,y_n,s_n)\\
        =& D_n \cdot \int_{\mathcal{Y}}\int_{\mathbb{R}_+}\mu_{n-1}(dy_{n-1},ds_{n-1}|h_{n-1})\\
        &\  \ \cdot \int_{\mathcal{Y}}\int_{\mathcal{X}}\nu(dy_n)\lambda(dx_n)q(x_n,y_n|x_{n-1},y_{n-1},g_{n-1}) \\
        & \  \ \cdot v(h_{n-1},g_{n-1},x_n,y_n,s_{n-1}+\beta^{n-1}c(x_{n-1},y_{n-1},g_{n-1})).
    \end{aligned}
\end{equation*}
The right-hand side of (\ref{eq:condition for info state recursion}) can be written as
\begin{equation*}
    \begin{aligned}
         \mathbb{E}_x^\pi[v'(h_{n-1},A_{n-1}, X_n)]
        = & \int_{\mathcal{X}}\lambda(dx_n)q^X(x_n|x_{n-1},y_{n-1},g_{n-1})\\
        & \  \ \cdot \int_{\mathcal{Y}}\int_{\mathbb{R}_+}\mu_{n-1}(dy_{n-1},ds_{n-1}|h_{n-1})v'(h_{n-1},g_{n-1},x_n)\\
        =&
        \int_{\mathcal{Y}}\mu_{n-1}^Y(dy_{n-1}|h_{n-1})\int_{\mathcal{X}}\lambda(dx_n)q^X(x_n|x_{n-1},y_{n-1},g_{n-1})\\
        & \  \ \cdot \int_{\mathcal{Y}}\int_{\mathbb{R}_+}v(h_{n-1},g_{n-1},x_n,y_n,s_n)\mu_n(dy_n,ds_n|h_n)\\
        =&
        \int_{\mathcal{Y}}\mu_{n-1}^Y(dy_{n-1}|h_{n-1})\int_{\mathcal{X}}\lambda(dx_n)q^X(x_n|x_{n-1},y_{n-1},g_{n-1})\\
        & \  \ \cdot \int_{\mathcal{Y}}\int_{\mathbb{R}_+}v(h_{n-1},g_{n-1},x_n,y_n,s_n)\nu(dy_n)q(x_n,y_n|x_{n-1},y_{n-1},g_{n-1})\\
        & \  \ \cdot \int_{\mathcal{Y}}\int_{\mathbb{R}_+}\delta_{s_{n-1}+\beta^{n-1}c(x_{n-1},y_{n-1},g_{n-1})}(ds_n)\mu_{n-1}(dy_{n-1},ds_{n-1}|h_{n-1})
        \\
        =&
        D_n\cdot \int_{\mathcal{Y}}\int_{\mathbb{R}_+}\mu_{n-1}(dy_{n-1},ds_{n-1}|h_{n-1})\\
        & \  \ \cdot \int_{\mathcal{Y}}\int_{\mathcal{X}}\nu(dy_n)\lambda(dx_n)q(x_n,y_n|x_{n-1},y_{n-1},g_{n-1})\\
        & \  \ \cdot v(h_{n-1},g_{n-1},x_n,y_n,s_{n-1}+\beta^{n-1}c(x_{n-1},y_{n-1},g_{n-1})).
    \end{aligned}
\end{equation*}
Hence we arrive at (\ref{eq:condition for info state recursion}).
Letting $v=\1_{B\times C}$ with $B\in\mathcal{B(\mathcal{Y}\times\mathbb{R}_+)}$ and $C\subset \mathcal{X\times A\times \cdots \times \mathcal{X}}$ be a measurable set of observable histories, we arrive at
\begin{equation*}
    \prod_{k=0}^n D_n\cdot\mathbb{P}_x^\pi((Y_n,S_n)\in B, (X_0,A_0,\cdots, X_n)\in C)=\mathbb{E}_x^\pi[\mu_n(B|X_0,A_0,\cdots,X_n)\1_{C}((X_0,A_0,\cdots,X_n))].
\end{equation*}
Consequently, we conclude that, given observable history $(X_0,A_0,\cdots,X_n)$, $\mu_n(\cdot|X_0,A_0,\cdots, X_n)$ encodes the conditional $\mathbb{P}_x^\pi$-distribution of $(Y_n,S_n)$ up to $\prod_{k=0}^n D_k$.
\qed
}

\textcolor{black}{
As can be observed from (\ref{eq:info state operator}), the dynamics of the unnormalized conditional density is a linear functional of the initial density. 
This linear structure simplifies the construction of the equivalent MDP for solving the DM's optimization problem which we will introduce next and the analysis of the impact of information manipulation that will be discussed in Section \ref{sec:impacts}.
}

\subsection{Decision-Maker's Dynamic Programming Equations}
\label{sec:solution:DP equations}
\textcolor{black}{
Based on the fact that the information states obtained from (\ref{eq:info state recursion}) specify the unnormalized conditional density of the unobservable states and the accumulated costs, we now construct an MDP to solve the DM's problem.
}
Define for $x\in\mathcal{X}$, $\mu\in\mathcal{M}(\mathcal{Y},\mathbb{R}_+)$, $z\in(0,1]$, and $n=1,\cdots,N$:
\begin{equation}
    \begin{aligned}
        V_{n\pi}(x,\mu,z)&:=\int_{\mathcal{Y}}\int_{\mathbb{R}_+}\mathbb{E}^\pi_{xy}\left[ U(s+z\sum_{k=0}^{n-1}\beta^kc(X_k,Y_k,A_k)) \right] \mu(dy,ds), \\
        V_n(x,\mu,z)&:=\inf_{\pi}V_{n,\pi}(x,\mu,z).
        \label{eq:V func}
    \end{aligned}
\end{equation}
\textcolor{black}{
In particular, we have $V_N(x, Q_0^Y\times \delta_0,1)=J_N(x)$.
}

\textcolor{black}{
Consider an MPD with state space $\mathcal{S}:=\mathcal{X}\times \mathcal{M}(\mathcal{Y}\times \mathbb{R}_+)\times (0,1]$, action space $\mathcal{A}$, and admissible actions from set $\mathcal{D}.$
The transition law of the MDP $\Tilde{Q}(\cdot|x,\mu,z,a)$ is specified for $(x,\mu,z,a)\in \mathcal{S}\times \mathcal{A}$, $a\in\mathcal{D}(x)$, and a measurable set $B\subset \mathcal{S}$ as $\Tilde{Q}(B|x,\mu,z,a):=\int_{\mathcal{X}}\1_{B} \left( (x',\Psi(x,a,x',\mu,z),\beta z) \right)\lambda(dx')$.
Let the stage cost be $0$ and the terminal cost be $V_0(x,\mu,z):=\int_{\mathcal{Y}}\int_{\mathbb{R}_+}U(s)\mu(dy,ds)$.
Define the strategy $f_n:\mathcal{X}\rightarrow \mathcal{A}$ such that $f_n(x,\mu,z)\in \mathcal{D}(x)$. 
Consider $f_n\in F$ and $\pi=(f_0,f_1,\cdots)\in\Pi^M$, where $F$ denotes the set of strategies and $\Pi^M$ denotes the set of Markov policies. 
Define the set $\mathcal{V}(\mathcal{S}):=\{v:\mathcal{S}\rightarrow\mathbb{R}: v \text{ is lower semicontinuous and }v\geq V_0 \}$, where we adopt the topology of weak convergence on $\mathcal{M}(\mathcal{Y}\times\mathbb{R}_+)$.
For $v\in\mathcal{V}(\mathcal{S})$ and $f\in F$ define the operator 
\begin{equation*}
    (T_fv)(x,\mu,z):=\int_{\mathcal{X}}v(x',\Psi(x,f(x,\mu,z),x',\mu,z),\beta z)\lambda(dx'), \  \ (x,\mu,z)\in \mathcal{S}.
\end{equation*}
The following result shows the procedure of obtaining the optimal policy from the MDP defined above. 
\begin{theorem}
\label{thm:verification thm}
The following statements hold.
\\
(a) For $\pi=(f_0,f_1,f_2,\cdots)\in\Pi^M$, the cost iteration admits $V_{n\pi}=T_{f_0}\cdots T_{f_{n-1}}V_0$ for $n=1\cdots,N$.
\\
(b) $V_n\in\mathcal{V}(\mathcal{S})$ and for $n=1,\cdots,N$ we have
\begin{equation}
    \label{eq:V func recursion}
    V_{n+1}(x,\mu,z)=\inf_{a\in\mathcal{D}(x)}\int_{\mathcal{X}}V_n(x',\Psi(x,a,x',\mu,z),\beta z)\lambda(dx'), \  \ (x,\mu,z )\in\mathcal{S}.
\end{equation}
(c) For $n=1,\cdots,N$, there exists a minimizer $f_n^*\in F$ of $V_{n-1}$. The optimal policy $(g_0^*,\cdots, g_{N-1}^*)$ of $J_N(x)$ is given by $g_n^*(h_n):=f_{N-n}^*(x_n,\mu_n(\cdot|h_n),\beta^n)$ for $n=0,\cdots,N-1$.
\end{theorem}
}
\textcolor{black}{
\noindent{\bf Proof}. 
We address part (a) using induction.
For $n=1$, let $a:=f_0(x,\mu,z)$.
Then,
\begin{equation*}
    \begin{aligned}
        T_{f_0}V_0(x,\mu,z)=& \int_{\mathcal{X}}V_0(x',\Psi(x,a,x',\mu,z),\beta z)\lambda(dx')\\
        =& \int_{\mathcal{X}} \int_{\mathcal{Y}}\int_{\mathbb{R}_+} U(s') \Psi(x,a,x',\mu,z)(dy',ds')\lambda(dx')
        \\
        =& \int_{\mathcal{X}}\int_{\mathcal{Y}}\int_{\mathbb{R}_+}\int_{\mathcal{Y}}\int_{\mathbb{R}_+}U(s')q(x',y'|x,y,a)\nu(dy')\delta_{s+zc(x,y,a)}(ds')\mu(dy,ds)\lambda(dx')
        \\
        =& \int_{\mathcal{Y}}\int_{\mathbb{R}_+}U(s+zc(x,y,a))\mu(dy,ds)\\
        =& V_{1\pi}(x,\mu,z).
    \end{aligned}
\end{equation*}
Suppose that the statement is true for $V_{n\pi}$, \ie, $V_{n\pi}=T_{f_0}\cdots T_{f_{n-1}}V_0$.
Let $a:=f_0(x,\mu,z)$ and $\Vec{\pi}=(f_1,f_2,\cdots)$  denotes the shifted policy for $\pi=(f_0,f_1,f_2,\cdots)\in \Pi^M$.
Then, 
\begin{equation*}
    \begin{aligned}
        T_{f_0}\cdots T_{f_n}V_0(x,\mu,z)
        =& T_{f_0}\left(T_{f_1}\cdots T_{f_n}\right)V_0(x,\mu,z)\\
        =&\int_{\mathcal{X}}V_{n\Vec{\pi}}(x',\Psi(x,f_0,x',\mu,z),\beta z)\lambda(dx')\\
        =&\int_{\mathcal{X}}\int_{\mathcal{Y}}\int_{\mathbb{R}_+}\mathbb{E}_{x'y'}^{\Vec{\pi}}\left[ U(s'+z\sum_{k=0}^{n-1}\beta^{k+1}c(X_k,Y_k,A_k)) \right]\\
        \  \ &\cdot \Psi(x,f_0,x',\mu,z)(dy',ds')\lambda(dx')\\
        =& \int_{\mathcal{X}}\int_{\mathcal{Y}}\int_{\mathbb{R}_+}\mathbb{E}_{x'y'}^{\Vec{\pi}}\left[ U(s'+z\sum_{k=0}^{n-1}\beta^{k+1}c(X_k,Y_k,A_k)) \right]
        \\
        \  \ & \cdot \int_{\mathcal{Y}}\int_{\mathbb{R}_+}q(x',y'|x,y,a)\nu(dy')\delta_{s+zc(x,y,a)}(ds')\mu(dy,ds)\lambda(dx')\\
        =& \int_{\mathcal{X}}\int_{\mathcal{Y}}\int_{\mathbb{R}_+} \mathbb{E}^{\pi}\left[ U(s'+z\sum_{k=1}^{n}\beta^{k}c(X_k,Y_k,A_k)) | X_1=x', Y_1=y' \right] \\
        \  \ & \cdot \delta_{s+zc(x,y,a)}(ds')\int_{\mathcal{Y}}\int_{\mathbb{R}_+}q(x',y'|x,y,a)\lambda(dx')\nu(dy')\mu(dy,ds)\\
        =& \int_{\mathcal{X}}\int_{\mathcal{Y}}\int_{\mathcal{Y}}\int_{\mathbb{R}_+} \mathbb{E}^{\pi}\left[ U(s+zc(x,y,a)+z\sum_{k=1}^{n}\beta^{k}c(X_k,Y_k,A_k)) | X_1=x', Y_1=y' \right]\\
        \  \ & \cdot q(x',y'|x,y,a)\lambda(dx')\nu(dy')\mu(dy,ds)\\
        =& \int_{\mathcal{Y}}\int_{\mathbb{R}_+}\mathbb{E}_{xy}^{\pi}\left[ U(s+z\sum_{k=0}^{n}\beta^k c(X_k,Y_k,A_k)) \right]\mu(dy,ds)\\
        =& V_{n+1 \pi} (x,\mu,z).
    \end{aligned}
\end{equation*}
The proofs of part (b) and (c) follow similar arguments as stated in the proof of Theorem 2 in \cite{baauerle2017partially}. 
The distinction lies in that we can skip the continuity property of $Q^X(\cdot|x,\mu^Y,a)$, which is defined by 
\begin{equation*}
    Q^X(B|x,\mu^Y,a):=\int_{B}\int_{\mathcal{Y}}q^X(x'|x,y,a)\mu^Y(dy)\lambda(dx'), B\in\mathcal{B}(\mathcal{X}).
\end{equation*}
The reason lies in that $Q^X(x_n|x_{n-1},\mu_{n-1}^Y,a_{n-1})=D_n$ and it is absent in both the update rule $\Psi$ and the transition rule $\Tilde{Q}$.
\qed
}

\textcolor{black}{
Compared with the setting adopted in \cite{baauerle2017partially}, the unnormalized conditional densities obtained using recursion (\ref{eq:info state recursion}) help simplify the transition law $\Tilde{Q}$ in the constructed MDP above.
The idea of adopting unnormalized conditional densities as information states is motivated by the observation that, in the proof of the cost iteration in Theorem 2 of \cite{baauerle2017partially}, the term $D_n$ is a common factor which appears in both the denominator due to the normalization of information states and the numerator due to the transition law of the construct MDP.
}

\subsection{Solution to Information Manipulation}
\label{sec:solution:solve info manipulation}
\textcolor{black}{
Based on the optimal policy obtained from Section \ref{sec:solution:DP equations}, the IM designs the joint state distributions for manipulation.
Note that the IM can always reconstruct the DM's optimal policy, since the IM is aware of all information that is available to the DM.
In this section, we first focus on the IM's design problem in one stage.
Then, we synthesize the solution to the DIMG building on a dynamic programming argument.
}

\textcolor{black}{
Let us elaborate the set of stagewise consistent designs $\Theta$.
We define $\Bar{\mathcal{P}}_n(x_{n-1},y_{n-1},a_{n-1}):=\{P_n\in\mathcal{P(\mathcal{X\times\mathcal{Y}}}): P_n^Y(\cdot)=Q^Y(\cdot|x_{n-1},y_{n-1},a_{n-1})  \}$ for $n=1,\cdots,N-1$, which denotes the set of stagewise consistent designs for the IM given the realizations of states and action of the previous stage.
Let $\Bar{\mathcal{P}}_0:= \{P_0 \in \mathcal{P}(\mathcal{X}\times \mathcal{Y}): P_0^Y=Q_0^Y \}$ denote the set of consistent designs at stage $n=0$.
Then, $\Theta=\Bar{\mathcal{P}}_0\times\cdots\times \Bar{\mathcal{P}}_{N-1}$.
Let $\pi^*=(g_0^*,\cdots,g_{N-1}^*)$ denote an optimal policy of $J_N(x)$.
Let $W_{0\pi^*}(x,y,a)=0$ for all $(x,y,a)\in\mathcal{X}\times\mathcal{Y}\times\mathcal{A}$.
For $n=N-1,\cdots,1$, define the following equations for $(x_{n-1},y_{n-1},a_{n-1})\in \mathcal{X}\times \mathcal{Y}\times \mathcal{A}$ and $\pi^*$:
\begin{equation}
    \begin{aligned}
        W_{N-n\pi^*}(x_{n-1},y_{n-1},a_{n-1}):=\inf_{P_{n}\in\Bar{\mathcal{P}}_{n}(x_{n-1},y_{n-1},a_{n-1})}
       &\alpha^{n}\left(  \int_{\mathcal{X}}\int_{\mathcal{Y}}r(x_{n},y_{n},g_{n}^*) P_{n}(dx_{n},dy_{n}) + \rho_{n}(P_{n}) \right)\\
       &\  \ + \int_{\mathcal{X}}\int_{\mathcal{Y}} W_{N-n-1 \pi^*}(x_n,y_n,g_n^*) P_n(dx_n,dy_n),
       \label{eq: W_N-n n=1,...,N-1}
    \end{aligned}
\end{equation}
and 
\begin{equation}
\begin{aligned}
    W_{N\pi^*}:=\inf_{P_{0}\in\Bar{\mathcal{P}}_{0}}  
    &\int_{\mathcal{X}}\int_{\mathcal{Y}}r(x_{0},y_{0},g_{0}^*) P_{0}(dx_{0},dy_{0}) + \rho_{0}(P_{0}) + \int_{\mathcal{X}}\int_{\mathcal{Y}} W_{N-1 \pi^*}(x_0,y_0,g_0^*) P_0(dx_0,dy_0).
       \label{eq:W_N}
\end{aligned}
\end{equation}
In (\ref{eq: W_N-n n=1,...,N-1}), the optimal manipulation at each stage depends on the realizations of the states and action from the previous stage $(x_{n-1},y_{n-1},a_{n-1})$. 
This dependence is due to the fact that the current design has to be chosen from the set of consistent designs $\Bar{\mathcal{P}}_n(x_{n-1},y_{n-1},a_{n-1})$ and that the IM recalls the whole history of states and actions up to the current stage.
The fact that we write the equations (\ref{eq: W_N-n n=1,...,N-1}) and (\ref{eq:W_N}) backward in time is also a consequence of stagewise consistency.
This can be observed since the optimal cost-to-go $W_{N-n-1\pi^*}$ depends on the design for the next stage $P_{n+1}$, which again has to be selected from the set of consistent designs $\Bar{\mathcal{P}}_{n+1}(x_n,y_n,a_n)$.
Moreover, the optimal manipulation depends on the current optimal control strategy $g_n^*$ rather than its implemented control action $a_n=g_n^*((h_{n-1},x_n))$ for the following reason.
The IM seeks to stealthily influence the DM's observation at the current stage $x_n$, and hence the control action $a_n=g_n^*((h_{n-1},x_n))$, by designing the probabilistic rule $P_n(dx_n,dy_n)$.
This stealthiness is achieved because the DM's  policy $\pi^*$ is derived from the original system model $Q$, meaning the DM remains unaware of the manipulation.
}

\textcolor{black}{
Let us also introduce the counterpart of (\ref{eq: W_N-n n=1,...,N-1}) and (\ref{eq:W_N}) under the interim design scheme.
In this scheme, the IM is not able to design the whole joint distribution on the state pair, since the true unobservable state has been realized.
Accordingly, stagewise consistency collapses since the IM only designs the distribution over the observable state.
Nevertheless, the evolution of the unobservable state is still governed by the true transition kernel, which can be observed from the following set of equations.
Let $W_{0\pi^*}^Y(y)=0$ for all $y\in\mathcal{Y}$. 
For $n=N-1,\cdots,1,0$, define the following equations for $y\in  \mathcal{Y}$ and $\pi^*$:
\begin{equation}
    \begin{aligned}
        W_{N-n\pi^*}^Y(y):=\inf_{P_{n}^X\in\mathcal{P}(\mathcal{X})}
       &\alpha^{n}\left(  \int_{\mathcal{X}}\int_{\mathcal{Y}}r(x_{n},y,g_{n}^*) P_{n}^X(dx_{n}) + \rho_{n}^X(P_{n}^X) \right)\\
       &\  \ + \int_{\mathcal{X}}\int_{\mathcal{Y}} W_{N-n-1 \pi^*}^Y(y') Q^Y(dy'|x_{n},y,g_{n}^*)P_n^X(dx_{n}),
       \label{eq:interim DP equations}
    \end{aligned}
\end{equation}
where $\rho_n^X(P_n^X)=2d_{TV}(P_n^X, Q^X(\cdot|x_{n-1},y_{n-1},a_{n-1}))$.
Note that in (\ref{eq:interim DP equations}), the value $W_{N-n\pi}^Y$ under the interim scheme does not depend on $(x_{n-1},y_{n-1},a_{n-1})$ as in (\ref{eq: W_N-n n=1,...,N-1}) under the ex ante scheme, since the choice of $P_n^X$ in (\ref{eq:interim DP equations}) is not constrained by stagewise consistency.
Instead, since IM observes the state realizations, his decision depends explicitly on the current unobservable state $y$, which not only influences the cost at the current stage but also the future costs via the transition rule $Q$.
Hence, equations (\ref{eq:interim DP equations}) are also formulated backward in time.
In the sequel, we assume that the solution sets of (\ref{eq:interim DP equations}) are nonempty and compact for all values of the unobservable state.
}

In the following, we first present a technical lemma, which is built on Theorem 2.2 in \cite{lasserre2010joint+}, to pave the way for showing the existence of solutions to DIMG in the ex ante and the interim schemes.

\begin{lemma}
\label{lemma:extension of Lasserre's thm}
Let $l:\mathcal{X}\times\mathcal{Y}\rightarrow\mathbb{R}$ denote a bounded function.
Let $P,Q \in \mathcal{P}(\mathcal{X}\times\mathcal{Y})$ denote joint probability measures on $\mathcal{X}\times\mathcal{Y}$ having measurable densities $p$ and $q$ with respect to probability measures $\lambda\in\mathcal{P}(\mathcal{X})$ and $\nu\in\mathcal{P}(\mathcal{Y})$.
Let $m\in \mathcal{P}(Y)$ denote a probability measure on $Y$.
Let the stochastic kernels $\phi,\psi$ from $\mathcal{Y}$ to $\mathcal{X}$ be defined by the following disintegrations given marginal $m(\cdot)$:
\begin{equation*}
\begin{aligned}
    &P(dx,dy)=\phi(dx|y)m(dy),\\
    &Q(dx,dy)=\psi(dx|y)m(dy).
\end{aligned}
\end{equation*}
Consider the optimization problem
\begin{equation}
    \begin{aligned}
        R:= \inf_{P\in\mathcal{P}(\mathcal{X}\times\mathcal{Y})}
        &\biggl\{ \int_{\mathcal{X}\times\mathcal{Y}} ldP + \int_{\mathcal{X}\times\mathcal{Y}} |dP-dQ|
        : \\
        &\  \
        P(\mathcal{X},B)=Q(\mathcal{X},B)=m(B), \forall B\in\mathcal{B(\mathcal{Y})} \biggr\},
        \label{eq:joint measure opt prob}
    \end{aligned}
\end{equation}
and, for all $y\in\mathcal{Y}$, the optimization problem
\begin{equation}
    T(y):=\inf_{\eta\in\mathcal{P}(\mathcal{X})}\int_{\mathcal{X}}l(x,y)\eta(dx)+\int_{\mathcal{X}}|\eta-\psi(\cdot|y)|dx.
    \label{eq:parametric opt prob}
\end{equation}
Suppose that the solution set $\mathcal{P}^*_y\subset \mathcal{P}(\mathcal{X})$ of (\ref{eq:parametric opt prob}) is nonempty and compact for all $y\in\mathcal{Y}$. Then, (\ref{eq:joint measure opt prob}) has an optimal solution and $R=\int_{\mathcal{Y}}T(y)m(dy)$.
\end{lemma}
\noindent{\bf Proof}.   
Since $\mathcal{P}^*_y$ is nonempty for all $y\in\mathcal{Y}$, we have
\begin{equation}
    \int_{\mathcal{X}}l(x,y)\eta_y(dx)+\int_{\mathcal{X}}|\eta_y-\psi(\cdot|y)|dx\geq T(y), \forall \eta_y\in\mathcal{P}(\mathcal{X}).
    \label{eq:proof:lemma:extension:1}
\end{equation}
With $\eta_y$ for all $y\in\mathcal{Y}$, $\eta_y(dx)m(dy)\in\mathcal{P}(\mathcal{X}\times\mathcal{Y})$ denotes a joint probability measure on $\mathcal{X}\times\mathcal{Y}$ with $Y$-marginal $m$. 
Then, for all $P$ feasible for (\ref{eq:joint measure opt prob}), we can represent it using a properly chosen $\eta_y$ for all $y\in\mathcal{Y}$ and the marginal $m$.
From (\ref{eq:proof:lemma:extension:1}), we obtain $\forall \eta_y\in\mathcal{P}(\mathcal{X})$, $\forall y\in\mathcal{Y}$, that
\begin{equation*}
    \int_{\mathcal{X}\times\mathcal{Y}}l(x,y)\eta_y(dx)m(dy)+\int_{\mathcal{X}\times\mathcal{Y}}|\eta_ym-\psi(\cdot|y)m|dxdy\geq \int_{\mathcal{Y}}T(y)m(dy).
    \label{eq:proof:lemma:extension:2}
\end{equation*}
Observing that $Q(dx,dy)=\psi(dx|y)m(dy)$ and $\eta_y(dx)m(dy)$ can represent any joint probability distribution on $\mathcal{X}\times\mathcal{Y}$ with $Y$-marginal $m$, we obtain $R\geq \int_{\mathcal{Y}}T(y)m(dy)$.
\par
Since the set-valued mapping $\mathcal{P}^*_y:\mathcal{Y}\rightrightarrows \mathcal{P}(\mathcal{X})$ is compact-valued, we observe from Proposition 4.4 in \cite{lasserre2010joint+} that there exists a measurable selector $\zeta:\mathcal{Y}\rightarrow \mathcal{P}^*_y$, such that the following holds for all $y\in\mathcal{Y}$,
\begin{equation*}
    \int_{\mathcal{X}}l(x,y)[\zeta(y)](dx)+\int_{\mathcal{X}}|[\zeta(y)](\cdot)-\psi(\cdot|y)|dx=T(y).
\end{equation*}
Pick $\mathcal{P}(\mathcal{X}\times\mathcal{Y})\ni P(dx,dy):=[\zeta(y)](dx)m(dy)$.
Since $P$ is feasible for (\ref{eq:joint measure opt prob}), the following holds
\begin{equation*}
    \begin{aligned}
        R\leq &\int_{\mathcal{X}\times\mathcal{Y}}ldP+\int_{\mathcal{X}\times\mathcal{Y}}|dP-dQ|\\
        \leq & \int_{\mathcal{Y}}\left( \int_{\mathcal{X}}l(x,y)[\zeta(y)](dx)+\int_{\mathcal{X}}|[\zeta(y)](\cdot)-\psi(\cdot|y)|dx \right)m(dy) \\
        \leq& \int_{\mathcal{Y}}T(y)m(dy).
    \end{aligned}
\end{equation*}
Hence, we conclude that $R=\int_{\mathcal{Y}}T(y)m(dy)$.
\qed

The following result concerns the solution to DIMG in the ex ante scheme.
 \textcolor{black}{
\begin{theorem}
\label{thm:existence ex ante}
Let $\pi^*=(g_0^*,\cdots, g_{N-1}^*)$ with $g_n^*(h_n)=f_{N-n}^*(x_n,\mu_n(\cdot|h_n),\beta^n)$ denote an optimal policy of $J_N(x)$ obtained from (\ref{eq:V func recursion}).
There exist minimizers $P_0^*,\cdots, P_{N-1}^*$ of (\ref{eq: W_N-n n=1,...,N-1}) and (\ref{eq:W_N}).
Moreover, $\theta^*=(P_0^*,\cdots, P_{N-1}^*)$ and $\pi^*=(g_0^*,\cdots, g_{N-1}^*)$ constitute an ex ante equilibrium of information distortion in the DIMG involving the DM and the IM.
\end{theorem}
\noindent{\bf Proof}. 
We first prove the existence result by relating the equations (\ref{eq: W_N-n n=1,...,N-1}) with the equations (\ref{eq:interim DP equations}). 
In particular, we show via an induction argument  that $W_{N-n\pi^*}$ is in the form of (\ref{eq:joint measure opt prob}) and $W_{N-n\pi^*}^Y$ is in the form of (\ref{eq:parametric opt prob}) if we choose $l_n(x_n,y_n)=\alpha^n r(x_n,y_n,g_n^*)+W_{N-n-1\pi^*}(x_n,y_n,g_n^*)$ in Lemma \ref{lemma:extension of Lasserre's thm}.
The statement is straightforward for $n=N-1$, \ie, $W_{1\pi^*}$ and $W_{1\pi^*}^Y$ can be related with $l_{N-1}(x_{N-1},y_{N-1})=\alpha^{N-1}r(x_{N-1},y_{N-1},g_{N-1}^*)$. 
And Lemma \ref{lemma:extension of Lasserre's thm} indicates that $W_{1\pi^*}=\int W_{1\pi^*}^Y (y)Q^Y(dy|x_{N-2},y_{N-2},g_{N-2}^*)$.
Suppose that the statement is true for $n=n+1$, \ie, we choose $l_{n+1}(x_{n+1},y_{n+1})=\alpha^{n+1}r(x_{n+1},y_{n+1},g_{n+1}^*)+W_{N-n-2\pi^*}(x_{n+1},y_{n+1},g_{n+1}^*)$ and arrive at $W_{N-n-1\pi^*}=\int W_{N-n-1\pi^*}^Y(y_{n+1})Q^Y(dy_{n+1}|x_n,y_n,g_n^*)$. 
Now, for $n=n$ we have by (\ref{eq:interim DP equations}) and the condition obtained from step $n=n+1$ that
\begin{equation*}
    \begin{aligned}
        W_{N-n\pi^*}^Y(y_n)=
        &\inf_{P_n^X\in\mathcal{P}(\mathcal{X})} 
        \alpha^n \Big( \int_{\mathcal{X}} r(x_n,y_n,g_n^*)P_n^X(dx_n) + \rho_n^X(P_n^X) \Big) \\
        & \  \ +\int_{\mathcal{X}}\int_{\mathcal{Y}}W_{N-n-1\pi^*}^Y(y_{N+1})Q^Y(dy_{n+1}|x_n,y_n,g_n^*)P_n(dx_n)\\
        =&\inf_{P_n^X\in\mathcal{P}(\mathcal{X})} 
        \alpha^n \Big( \int_{\mathcal{X}} r(x_n,y_n,g_n^*)P_n^X(dx_n) + \rho_n^X(P_n^X) \Big) \\
        & \ \ + \int_{\mathcal{X}}W_{N-n-1\pi^*}(x_n,y_n,g_n^*)P_n^X(dx_n).
    \end{aligned}
\end{equation*}
This certifies that $l_n(x,y)=\alpha^nr(x,y,g_n^*)+W_{N-n-1\pi^*}(x,y,g_n^*)$ is the choice for $W_{N-n\pi^*}$ and $W_{N-n\pi^*}^Y$.
Therefore, Lemma \ref{lemma:extension of Lasserre's thm} indicates that minimizers $P_0^*,\cdots,P_{N-1}^*$ exist for the optimization problems in (\ref{eq: W_N-n n=1,...,N-1}) and (\ref{eq:W_N}), and the optimal values satisfy
\begin{equation}
    W_{N-n\pi^*}=\int W_{N-n\pi^*}^Y Q^Y(dy_n|x_{n-1},y_{n-1},a_{n-1}), \  \ n=1,\cdots,N-1,
    \label{eq:relation optimal values ex ante and interim}
\end{equation}
and $W_{N\pi^*}=\int W_{N\pi^*}^Y Q^Y_0(dy_0)$.
\\
To show that $\pi^*$ and $\theta^*$ constitute an ex ante equilibrium of information distortion, we first elaborate the probability measure $\mathbb{P}^{\pi\theta}$.
Let $P_n^{X|Y}(\cdot|y)$ denote the conditional probability of the observable state given unobservable state $y$ associated with $P_n$.
Then, we obtain for $\pi=(g_0,\cdots,g_{N-1})$ and $\theta=(P_0,\cdots, P_{N-1})$ that
\begin{equation}
\begin{aligned}
    &\mathbb{P}^{\pi\theta}((dx_0,dy_0),\cdots,(dx_{N-1},dy_{N-1}))\\
    =& P_0(dx_0,dy_0)\cdots P_{N-1}(dx_{N-1},dy_{N-1})\\
    =& P_0^{X|Y}(dx_0|y_0)P_0^Y(dy_0)\cdots P_{N-1}^{X|Y}(dx_{n-1}|y_{N-1})P_{N-1}^Y(dy_{N-1})\\
    =& P_0^{X|Y}(dx_0|y_0)Q^Y(dy_0)\cdots P_{N-1}^{X|Y}(dx_{n-1}|y_{N-1})Q^Y(dy_{N-1}|x_{N-2},y_{N-2},g_{N-2}),
    \label{eq:proof:P pi theta measure}
\end{aligned}
\end{equation}
where the last equality is a consequence of stagewise consistency.
Based on (\ref{eq:proof:P pi theta measure}), we can express the IM's optimization problem as the following nested form:
\begin{equation}
\begin{aligned}
    I_{N\pi}=&
    \inf_{\theta\in\Theta}\int\cdots\int\left[ \sum_{k=0}^{N-1}\alpha^k \left( r(X_k,Y_k,A_k)+ \rho_k(P_k) \right) \right]\mathbb{P}^{\pi\theta}((dx_0,dy_0),\cdots,(dx_{N-1},dy_{N-1}))\\
    =& \inf_{\theta\in\Theta} \int\Bigg( \alpha^0\Big(r(x_0,y_0,g_0)+\rho_0(P_0)\Big) + \int \Bigg( \alpha^1\Big( r(x_1,y_1,g_1)+\rho_1(P_1) \Big)+\cdots
    \\
    & \  \ + 
    \int \Bigg( \alpha^{N-1}\Big(r(x_{N-1},y_{N-1},g_{N-1}) +\rho_{N-1}(P_{N-1})\Big)
    \Bigg) P_{N-1}(dx_{N-1},dy_{N-1}) \\
    & \  \ \cdots \Bigg) P_1(dx_1,dy_1)\Bigg)P_0(dx_0,dy_0)\\
    =& \inf_{P_0\in\Bar{\mathcal{P}}_0} \int\Bigg( \alpha^0\Big(r(x_0,y_0,g_0)+\rho_0(P_0)\Big) + \inf_{P_1\in\Bar{\mathcal{P}}_1}\int \Bigg( \alpha^1\Big( r(x_1,y_1,g_1)+\rho_1(P_1) \Big)+\cdots
    \\
    & \  \ + \inf_{P_{N-1}\in\Bar{\mathcal{P}}_{N-1}}
    \int \Bigg( \alpha^{N-1}\Big(r(x_{N-1},y_{N-1},g_{N-1}) +\rho_{N-1}(P_{N-1})\Big)
    \Bigg) P_{N-1}(dx_{N-1},dy_{N-1}) \\
    & \  \ \cdots \Bigg) P_1(dx_1,dy_1)\Bigg)P_0(dx_0,dy_0).
    \label{eq:proof:I_npi nested form}
\end{aligned}
\end{equation}
The last equation in (\ref{eq:proof:I_npi nested form}) and (\ref{eq: W_N-n n=1,...,N-1}) lead to the relation $I_{N\pi^*}=W_{N\pi^*}$.
Hence, $\theta^*=(P_0^*,\cdots,P_{N-1}^*)$ obtained by solving (\ref{eq: W_N-n n=1,...,N-1}) and (\ref{eq:W_N}) satisfies $I_{N\pi^*}=I_{N\pi^*\theta^*}$. 
Therefore, $\pi^*$ and $\theta^*$ constitutes and ex ante equilibrium of information distortion.
\qed
}

\textcolor{black}{
Note that the set of equations in (\ref{eq: W_N-n n=1,...,N-1}) and (\ref{eq:W_N}) are not exactly the same as the dynamic programming equations for optimal control problems.
The major difference lies in that, in (\ref{eq: W_N-n n=1,...,N-1}) and (\ref{eq:W_N}), the IM does not observe the state realizations prior to making a decision at each stage.
Instead, the IM designs ex ante the joint distribution of the states before their realizations at each stage. 
One also observes that equations (\ref{eq: W_N-n n=1,...,N-1}) and (\ref{eq:W_N}) are formulated and solved according to backward induction, which is similar as in dynamic programming equations.
This is because the designs of the IM needs to take into account not only the current cost but also the cost-to-go for the subsequent stages, which appear in (\ref{eq: W_N-n n=1,...,N-1}) and (\ref{eq:W_N}) due to stagewise consistency.
}

In the following, we connect the ex ante and the interim design schemes.
From the proof of Lemma \ref{lemma:extension of Lasserre's thm}, we observe that there exists a stochastic kernel $\phi^*_n$ from $\mathcal{Y}$ to $\mathcal{X}$ for every $P_n^*$ that attains the minimum of of $W_{N-n\pi^*}$ such that the following disintegration holds:
\begin{equation}
    P^*_n(B_x\times B_y)=\int_{B_y}\phi^*_n(B_x|y_n)Q^Y(dy_n|x_{n-1},y_{n-1},a_{n-1}), \  \ \forall B_x\in\mathcal{B}(\mathcal{X}), B_y\in\mathcal{B}(\mathcal{Y}).
    \label{eq:optimal design disintegration}
\end{equation}
Furthermore, by (\ref{eq:relation optimal values ex ante and interim}), we arrive at the following result that is a direct consequence of the existence of ex ante equilibrium of information distortion and (\ref{eq:optimal design disintegration}).
\begin{corollary}
\label{coro:existence stackelberg interim}
The sequences $(\phi^*_0,\phi^*_1,\cdots,\phi^*_{N-1})$ and $(g_0^*,g_1^*,\cdots,g_{N-1}^*)$ constitute an interim equilibrium of information distortion in the DIMG.
\end{corollary}

Note that the interim solutions $(\phi^*_0,\phi^*_1,\cdots,\phi^*_{N-1})$ can be regarded as degenerate cases of the ex ante solutions $(P_0^*,P_1^*,\cdots,P_{N-1}^*)$.
However, both manipulation schemes possess distinct advantages.
The ex ante solution can be fully deployed once the IM obtains the optimal policy of the DM $\pi^*$.
The manipulated kernels $P_0^*, \cdots, P_{N-1}^*$ then automatically generate the state pair at each stage.
Since the IM does not have to obtain any output from the true kernel in this scheme, he possesses a higher level of stealthiness.
The interim solution requires the unobservable states as input for manipulation, which has the interpretation that the IM trades stealthiness for manipulation effects.

\subsection{Examples}
\label{sec:example}
In this section, we present two examples of the DIMG in the ex ante information pattern.
Both of them adopt two-stage ($N=2$) settings and we consider the scenario where the IM only performs information manipulation at stage $n=1$. 
The first example builds on discrete state spaces and aims to illustrate the procedure of information manipulation as well as tractable reformulations of the IM's stage problem.
The second example adopts a linear system with Gaussian noises. 
We provide insights for the optimal policies of the DM and the IM.
Cost functions, when there is a need to specify them, are assumed quadratic in both examples.
We also adopt the setting where $U(z)=z$ for notational convenience.

\paragraph{Discrete state space}
We restrict the unobservable state space to $\mathcal{Y}=\{y^1,y^2\}$.
We use $y_n^i$, for $i=1,2$, to denote the case where the realization of the unobservable state at stage $n$ is $y^i$.
\textcolor{black}{
The transition probability of the unobservable state, for $n=0,1$, and $0<a_n<1$, is defined by
\begin{equation}
    \varphi_y(y_{n+1}|y_n,a_n)=
    \begin{cases}
        1-a_n, \  \ &\text{if } y_{n+1}=y_n, \\
        a_n,\  \  &\text{otherwise}.
    \end{cases}
    \label{eq:eg:varphi y}
\end{equation}
From (\ref{eq:eg:varphi y}) we observe that the action represents the probability that the unobservable state switches.
}
The observable state $x$ is generated by the relation $x_{n+1}=y_n+\omega_{n+1}$ for $n=0,1$, where the noise process $\omega$ is independent and identically distributed for $n=0,1$  and has two possible realization $\omega^1$ and $\omega^2$ with probabilities $\tau\in (0,1)$ and $1-\tau$, respectively.
The two possible values of the unobservable state and the two possible values of the noise generate four potential candidates of the observable state, \ie, $x^{i,j}$ for $i,j\in \{1,2\}$, where $x^{i,j}$ correspond to the candidate generated by $y^i$ and $\omega^j$.
We further assume that the values of the four potential candidates are distinct.
\textcolor{black}{
Let $\varphi_x(x_{n+1}|y_n)$ denote the realization probability of $x_{n+1}$ given $y_n$ for $n=0,1$, which is defined by
\begin{equation*}
    \begin{aligned}
        \varphi_x(x_{n+1}|y_n)=\begin{cases}
            \tau, \  \ &\text{if } y_n=y_n^i \text{ and } x_{n+1}=x_{n+1}^{i,1}, \forall i=1,2,
            \\
            1-\tau, \  \ &\text{if } y_n=y_n^i \text{ and } x_{n+1}=x_{n+1}^{i,2}, \forall i=1,2,
            \\
            0, \  \ &\text{otherwise}.
        \end{cases}
    \end{aligned}
\end{equation*}
}
The above setting of the control system specifies the transition law 
\begin{equation*}
Q(x_{n+1},y_{n+1}|y_n,a_n)=\varphi_y(y_{n+1}|y_n,a_n)\varphi_x(x_{n+1}|y_n).
\end{equation*}
Let $Q_0^Y=(Q_0^{y^1}, Q_0^{y^2})$ denote the initial distribution of the unobservable state.
The information state recursion (\ref{eq:info state recursion}), for $B\in\mathcal{B}(\mathbb{R})$, admits the following form 
\begin{equation*}
\begin{aligned}
    \mu_1(y_1 \times B | h_0,a_0,x_1)= &\int_{B}\varphi_y(y_1|y_0^1,a_0)\varphi_x(x_{1}|y_0^1)\delta_{s_0+c(y_0^1,a_0)}(ds_1)Q_0^{y^1}\delta_0(ds_0) \\
    \  \ &+
    \int_{B}\varphi_y(y_1|y_0^2,a_0)\varphi_x(x_{1}|y_0^2)\delta_{s_0+c(y_0^2,a_0)}(ds_1)Q_0^{y^2}\delta_0(ds_0),
\end{aligned}
\end{equation*}
\begin{equation*}
\begin{aligned}
    \mu_2(y_2 \times B | h_1,a_1,x_2)= &\int_{B}\varphi_y(y_2|y_1^1,a_1)\varphi_x(x_{2}|y_1^1)\delta_{s_1+c(y_1^1,a_1)}(ds_2)\mu_1(y_1^1, ds_1) \\
    \  \ &+
    \int_{B}\varphi_y(y_2|y_1^2,a_1)\varphi_x(x_{2}|y_1^2)\delta_{s_1+c(y_1^2,a_1)}(ds_2)\mu_1(y_1^2, ds_1). 
\end{aligned}
\end{equation*}

We choose the cost function of the DM as $c(y,a)=y^2+a^2+\hat{c}ay$ with $\hat{c}>0$.
According to (\ref{eq: W_N-n n=1,...,N-1}) and Theorem \ref{thm:existence ex ante}, we proceed backwards in time  starting with deriving the DM's optimal policy at stage $n=1$.
This is followed by determining the IM's optimal design at stage $n=1$, which relies on the DM's optimal policy.
By (\ref{eq:V func recursion}), we determine the optimal policy $a_1^*$ of stage $n=1$ by solving
\begin{equation*}
\begin{aligned}
    V_1(x_1,\mu_1,\beta) =& \inf_{a_1}\int_{\mathcal{X}}V_0(x_2,\mu_2,\beta^2)dx_2.
\end{aligned}
\end{equation*}
By elaborating $V_0(x_2,\mu_2,\beta^2)=\int\int U(s_2)\mu_2(dy_2,ds_2)$ and ignoring the terms that are independent of $a_1$, we obtain the following optimization problem for finding the DM's optimal policy:
\begin{equation}
\begin{aligned}
    a_1^*=\argmin_{a_1} \  \
    \bigg\{ &c(y_1^1,a_1)\left( \varphi_y(y_1^1|y_0^1,a_0)\varphi_x(x_{1}|y_0^1)Q_0^{y^1} +\varphi_y(y_1^1|y_0^2,a_0)\varphi_x(x_{1}|y_0^2)Q_0^{y^2} \right) \\
    & + c(y_1^2,a_1)\left( \varphi_y(y_1^2|y_0^1,a_0)\varphi_x(x_{1}|y_0^1)Q_0^{y^1} +\varphi_y(y_1^2|y_0^2,a_0)\varphi_x(x_{1}|y_0^2)Q_0^{y^2} \right).
    \label{eq:eg:a_1^* 1}
\end{aligned}
\end{equation}
To arrive at a compact representation, we define the following quantities:
\begin{equation*}
    \begin{aligned}
        &\mathfrak{m}^1(x_1,a_0):=\varphi_x(x_{1}|y_0^1)(1-a_0)Q_0^{y^1}+\varphi_x(x_{1}|y_0^2)a_0Q_0^{y^2},
        \\      &\mathfrak{m}^2(x_1,a_0):=\varphi_x(x_{1}|y_0^1)a_0Q_0^{y^1}+\varphi_x(x_{1}|y_0^2)(1-a_0)Q_0^{y^2}.
    \end{aligned}
\end{equation*}
Combining (\ref{eq:eg:varphi y}) and (\ref{eq:eg:a_1^* 1}), we obtain the DM's optimal policy 
\begin{equation}
    a_1^*(x_1,a_0)=\frac{-\hat{c}(y^1\mathfrak{m}^1(x_1,a_0)+y^2\mathfrak{m}^2(x_1,a_0))}{2(\mathfrak{m}^1(x_1,a_0)+\mathfrak{m}^2(x_1,a_0))},
    \label{eq:eg:a_1^* 2}
\end{equation}
which is a function of $x_1$ and $a_0$.

Observing $a_1^*$, the IM determines the joint distribution of $x_1$ and $y_1$ by resorting to (\ref{eq: W_N-n n=1,...,N-1}), which admits the following form:  
\begin{equation}
\begin{aligned}
    W_{1\pi\theta}&=\inf_{P_1\in\Bar{\mathcal{P}}(x_0,y_0,a_0)}\int_{\mathcal{X}}\int_{\mathcal{Y}}r(x_1,y_1,a_1^*(x_1,a_0))p_1(x_1,y_1)dx_1dy_1 + \rho_1(P_1).
    \label{eq:eg:IM prob at stage 1}
\end{aligned}
\end{equation}
In general, the IM's cost function $r$ depends on $x_1$, $y_1$, and $a_1^*$ simultaneously.
Since the optimal policy $a_1^*$ is a function of $x_1$, the dependence on $x_1$ of the IM's cost $r$ can be implicit through $a_1^*$.
When the cost $r$ is not separable in $x_1$ and $y_1$, the joint distribution $P_1$ contains more information than its 
marginals.
In this scenario, there is more space for the IM to perform ex ante information manipulation to steer the states and the DM's controls.

In the following, we assume that $r$ is not separable in $x_1$ and $y_1$, and further elaborate the IM's manipulation.
We specify the joint probability density density $p_1(x_1,y_1)$ using the following matrix:
\begin{equation*}
    p_1(x_1,y_1)=\begin{pmatrix}
        p_1^{1,1} \  \  p_1^{1,2} \ 
 \  p_1^{2,1} \  \  p_1^{2,2}
        \\
        p_2^{1,1} \  \  p_2^{1,2} \ 
 \ p_2^{2,1} \  \  p_2^{2,2}
    \end{pmatrix},
\end{equation*}
where $p_k^{i,j}$ denotes the probability of $(x_1^{i,j}, y_1^k)$ for $i,j,k=1,2$ and $\sum_{i,j,k=1,2}p_k^{i,j}=1$.
Accordingly, the expected cost of the IM in the objective function in (\ref{eq:eg:IM prob at stage 1}) can be expressed as
\begin{equation}
\int_{\mathcal{X}}\int_{\mathcal{Y}}r(x_1,y_1,a_1^*)p_1(x_1,y_1)dx_1dy_1 = \sum_{i=1,2}\sum_{j=1,2}\sum_{k=1,2}p_k^{i,j} r(x_1^{i,j},y_1^k, a_1^*(x_1^{i,j},a_0)).
\label{eq:eg:average cost terms}
\end{equation}

From (\ref{eq:eg:IM prob at stage 1}), we observe that the IM's optimal design depends on the state realizations $x_0$ and $y_0$ and the DM's control $a_0$ at the previous stage.
\textcolor{black}{
This setting is feasible since we have assumed that the IM has perfect observation about the states.
Then, for given $y_0,a_0$, the stagewise consistency constraint in (\ref{eq:eg:IM prob at stage 1}) reads 
\begin{equation*}
\begin{aligned}
    \begin{cases}
        p_1^{1,1}+p_1^{1,2}+p_1^{2,1}+p_1^{2,2}=\varphi_y(y_1^1|y_0,a_0),\\
        p_2^{1,1}+p_2^{1,2}+p_2^{2,1}+p_2^{2,2}=\varphi_y(y_1^2|y_0,a_0).
    \end{cases}
\end{aligned}
\end{equation*}
The manipulation cost in the objective function in (\ref{eq:eg:IM prob at stage 1}) becomes
\begin{equation}
\begin{aligned}
    \rho_1(P_1)=&
    |p_1^{1,1}-\varphi_y(y_1^1|y_0,a_0)\tau|
    +|p_1^{1,2}-\varphi_y(y_1^1|y_0,a_0)(1-\tau)|\\
    & \  \ 
    +|p_2^{2,1}-\varphi_y(y_1^2|y_0,a_0)(1-\tau)|
    +|p_2^{2,2}-\varphi_y(y_1^2|y_0,a_0)\tau|.
    \label{eq:eg:abolute values}
\end{aligned}
\end{equation}
The absolute value terms in (\ref{eq:eg:abolute values}) can be further replaced by introducing auxiliary variables $\mathfrak{a}_k^j$ for  $k,j=1,2$, to the optimization problem and considering corresponding auxiliary constraints.
Therefore, under the settings specified above, we can reformulate (\ref{eq:eg:IM prob at stage 1}) into the following linear programming problem:
\begin{equation*}
    \begin{aligned}
        \min_{p_1, \mathfrak{a}_k^j, k,j=1,2}
        \  \ &
         \sum_{i=1,2}\sum_{j=1,2}\sum_{k=1,2}p_k^{i,j} r(x_1^{i,j},y_1^k, a_1^*(x_1^{i,j},a_0))
        + \mathfrak{a}_1^1+ \mathfrak{a}_1^2+ \mathfrak{a}_2^1+ \mathfrak{a}_2^2
        \\
        \text{s.t.} \  \ &
        \sum_{i,j,k=1,2}p_k^{i,j}=1, \forall i,j,k=1,2, \\
        & p_1^{1,1}+p_1^{1,2}+p_1^{2,1}+p_1^{2,2}=\varphi_y(y_1^1|y_0,a_0), \\
        & p_2^{1,1}+p_2^{1,2}+p_2^{2,1}+p_2^{2,2}=\varphi_y(y_1^2|y_0,a_0),\\
        & 
        \mathfrak{a}_1^1\geq p_1^{1,1}-\varphi_y(y_1^1|y_0,a_0)\tau, \  \ \mathfrak{a}_1^1\geq -(p_1^{1,1}-\varphi_y(y_1^1|y_0,a_0)\tau), \\
        & \mathfrak{a}_1^2\geq p_1^{1,2}-\varphi_y(y_1^1|y_0,a_0)(1-\tau), \  \ \mathfrak{a}_1^2\geq -(p_1^{1,2}-\varphi_y(y_1^1|y_0,a_0)(1-\tau)),  \\
        & \mathfrak{a}_2^1\geq p_2^{2,1}-\varphi_y(y_1^2|y_0,a_0)(1-\tau), \  \ \mathfrak{a}_2^1\geq -(p_2^{2,1}-\varphi_y(y_1^2|y_0,a_0)(1-\tau)), \\
        &  \mathfrak{a}_2^2\geq p_2^{2,2}-\varphi_y(y_1^2|y_0,a_0)\tau, \  \ \mathfrak{a}_2^2\geq -(p_2^{2,2}-\varphi_y(y_1^2|y_0,a_0)\tau).
    \end{aligned}
\end{equation*}
}

\paragraph{Continuous state space}
Consider $\mathcal{X}=\mathcal{Y}=\mathbb{R}$ and the following linear system
\begin{equation}
    \begin{cases}
        x_{n+1} = hy_n + \omega_{n+1}, \\
        y_{n+1} = \Tilde{b}y_n+\hat{b}a_n + \epsilon_{n+1},
    \end{cases}
    \label{eq:eg:linear sys}
\end{equation}
where the process $x$ is observable and the process $y$ is unobservable, and $h,\Tilde{b},\hat{b}\in\mathbb{R}$, are system parameters.
The noise processes $\omega$ and $\epsilon$ are both independent and identically distributed standard normal distributions with probability density functions $\varphi_\omega$ and $\varphi_\epsilon$, respectively.
With the control system specified, we can write the transition law for $B_1,B_2 \in \mathcal{B}(\mathbb{R})$ as
\begin{equation*}
    Q(B_1\times B_2|x,y,a)=\int_{B_1}\varphi_{\omega}(x'-hy)dx'\int_{B_2}\varphi_{\epsilon}(y'-\Tilde{b}y-\hat{b}a)dy'.
\end{equation*}
Accordingly, the information state recursion (\ref{eq:info state recursion}), for $B_1,B_2 \in \mathcal{B}(\mathbb{R})$, admits the following form
\begin{equation*}
\begin{aligned}
    \mu_0(B_1\times B_2|h_0) = &Q_0^Y(B_1)\times \delta_0(B_2), \\
    \mu_1(B_1\times B_2|h_0,a_0,x_1) =& \int_{\mathcal{Y}}\int_{\mathbb{R}_+}\int_{B_1}\int_{B_2}\varphi_\omega(x_1-hy_0)\varphi_\epsilon(y_1-\Tilde{b}y_0-\hat{b}a_0)dy_1\\
    &\cdot\delta_{s_0+c(x_0,y_0,a_0)}(ds_1)\mu_0(dy_0, ds_0),\\
    \mu_2(B_1\times B_2|h_1,a_1,x_2) =& \int_{\mathcal{Y}}\int_{\mathbb{R}_+}\int_{B_1}\int_{B_2}\varphi_\omega(x_2-hy_1)\varphi_\epsilon(y_2-\Tilde{b}y_1-\hat{b}a_1)dy_2\\
    &\cdot\delta_{s_1+c(x_1,y_1,a_1)}(ds_2)\mu_1(dy_1, ds_1).
\end{aligned}
\end{equation*}
The cost functions are chosen as $c(y,a)=y^2+a^2+\hat{c}ay$ with $\hat{c}>0$ and $r(y,a)=y^2+a^2+\hat{r}ay$ with $\hat{r}>0$.
For notational simplicity, we adopt the setting where $Q_0^Y(y_0)$ is the standard normal distribution.
We also restrict the IM's action set to be in the class of jointly Gaussian distributions and let $\gamma=0$.

Following the same procedure shown in the discrete state space scenario, we determine the optimal policy $a_1^*$ of stage $n=1$ by resorting to (\ref{eq:V func recursion}), which reduces to
\begin{equation*}
\begin{aligned}
    V_1(x_1,\mu_1,\beta) =& \inf_{a_1}\int_{\mathcal{X}}V_0(x_2,\mu_2,\beta^2)dx_2 \\
    = &\inf_{a_1}\int_{\mathcal{Y}}\int_{\mathcal{Y}}\left[ c(y_0,a_0)+\beta c(y_1,a_1) \right] \varphi_{\omega}(x_1-hy_0)\\
    &\cdot\varphi_\epsilon(y_1-\Tilde{b}y_0-\hat{b}a_0)dy_1Q_0^Y(dy_0).
\end{aligned}
\end{equation*}
Let $\iota=\frac{\hat{c}\Tilde{b}h}{2(h^2+1)}$.
Using first-order optimality conditions and properties of Gaussian distributions, we arrive at $a_1^*=-\frac{\hat{c}\hat{b}}{2}a_0-\iota x_1$.
The maneuverability of the DM's optimal policy $a_1^*$ heavily depends on the system parameter $h$.
When $|h|\rightarrow +\infty$, one observes $|\iota|\rightarrow 0$.
In this scenario, the impact of the observation $x_1$ on the control specified by $a_1^*$ is infinitesimal.
Hence, it is more challenging to perform information manipulation.
From the system equations (\ref{eq:eg:linear sys}), we observe that a relatively small $|h|$ indicates that the noise $\omega$ plays a more dominant role in the observation kernel.
Therefore, a system is more susceptible to information manipulation if its observations are noisier.

Since we have restricted $p_1$ to be a jointly Gaussian distribution and its $Y$-marginal needs to admit density  $p_1^Y(y_1)=\varphi_\epsilon(y_1-\Tilde{b}y_0-\hat{b}a_0)$ according to the stage-wise consistency constraint, we identify $p_1$ with its $X$-marginal $p_1^X(x_1)=\frac{1}{\sqrt{2\pi\mathfrak{v}_1^2}}\text{exp}\left(\frac{(x_1-\mathfrak{m_1})^2}{-2\mathfrak{v}_1^2}\right)$ with parameters $\mathfrak{m_1}\in\mathbb{R}$ and $\mathfrak{v}_1>0$ and its correlation coefficient $|\mathfrak{r}_1|< 1$.
Then,  we can reformulate the IM's stage problem into the following optimization problem:
\begin{equation}
    \inf_{\mathfrak{m}_1, \mathfrak{v}_1,\mathfrak{r}_1}\  \ \iota^2(\mathfrak{m}_1^2+\mathfrak{v}_1^2)  + 
    (\hat{c}\hat{b}a_0 - \hat{r}(\Tilde{b}y_0+\hat{b}a_0))\iota \mathfrak{m}_1
    + \hat{r}\iota \mathfrak{r}_1\mathfrak{v}_1.
    \label{eq:example:IM prob at stage 1 Gaussian dist}
\end{equation}
An inspection on (\ref{eq:example:IM prob at stage 1 Gaussian dist}) leads to the fact that $\mathfrak{r}$ is chosen as either $1$ or $-1$ depending on the sign of $\hat{r}\iota\mathfrak{v}_1$.
Suppose that $\hat{r}\iota\mathfrak{v}_1>0$ and one chooses the optimal $\mathfrak{r}^*=-1$.
Then, the optimal solution of (\ref{eq:example:IM prob at stage 1 Gaussian dist}) admits $\mathfrak{v}_1^*=\frac{\hat{r}}{2\iota}$ and $\mathfrak{m}_1^*=\frac{\hat{c}\hat{b}a_0 - \hat{r}(\Tilde{b}y_0+\hat{b}a_0)}{-2\iota}$.
The coefficient of variation (CV) of the $X$-marginal of the optimal information design can be obtained as $\frac{\mathfrak{v}_1^*}{|\mathfrak{m}_1^*|}=\frac{\hat{r}}{|\hat{r}(\Tilde{b}y_0+\hat{b}a_0)-\hat{c}\hat{b}a_0|}$, which is a fixed quantity at stage $n=1$.
This fact indicates that the IM can optimally perform manipulation based on practical concerns, such as information generation difficulties and communication limitations, as long as the statistical dispersion of the designed information is properly controlled.
Furthermore, we observe that the CV is a decreasing function of $|\hat{r}-\hat{c}|$, which measures the degree of misalignment between the objectives of the IM and the DM.
This monotonicity indicates that the distribution of the optimally designed information is more concentrated when there is a significant misalignment between the objectives of the IM and the DM.
When the DM can observe multiple samples at one stage, or when the IM's design is similar or identical in successive stages, a more concentrated distribution results in a greater possibility for the DM to realize or detect from statistics of observations that the IM is present.

\section{Impacts of Information Manipulation}
\label{sec:impacts}
In this section, we investigate the impacts of information manipulation on the effectiveness of policy obtained by the DM.
In particular, we derive an upper-bound of the performance degradation of the DM's control efforts by comparing the scenario involving the IM with the scenario where the IM is absent.

\textcolor{black}{
In the sequel, we assume that the initial observation $x$ is drawn from a given distribution $Q_0^X$, which admits a density function $q_0^X$ with respect to $\lambda$.
Let $Q_0=Q_0^XQ_0^Y$.
We denote by $\mathbb{P}^\pi$ the counterpart of $\mathbb{P}^\pi_x$ when distributional initiations are considered.
We write $J_{N\pi}=\int J_{N\pi}(x)dQ_0^X(x)$.
}

\textcolor{black}{
Recall that the true information state process $\mu_0,\cdots,\mu_{N-1}$ should be generated by the operator (\ref{eq:info state operator}) if the DM receives observation samples from the true kernel $Q$.
When the IM is present, however, the DM receives observation samples from the manipulated kernels $P_0,\cdots,P_{N-1}$. 
We use $\Tilde{\mu}_0,\cdots,\Tilde{\mu}_{N-1}$ to denote the manipulated information states  generated by the manipulated operators $\Tilde{\Psi}_n:\mathcal{X}\times \mathcal{A}\times\mathcal{X}\times\mathcal{M}(\mathcal{Y}\times\mathbb{R}_+)\mathbb{R}_+\rightarrow\mathcal{M}(\mathcal{Y}\times\mathbb{R}_+)$ defined for $n=1,\cdots,N-1$ by
\begin{equation*}
    \Tilde{\Psi}_n(x,a,x',\Tilde{\mu},z)(B):=\int_{\mathcal{Y}}\int_{\mathbb{R}_+}\left( \int_{B}p_n(x',y')\nu(dy')\delta_{s+zc(x,y,a)}(ds') \right) \Tilde{\mu}(dy,ds), \  \ B\in\mathcal{B}(\mathcal{Y}\times\mathbb{R}_+).
    \label{eq:manipulated info state operator}
\end{equation*}
We are interested in comparing the performances of the DM's control policy between the scenario where the IM is present, \ie, when the DM adopts the manipulated information states $\Tilde{\mu}_0,\cdots,\Tilde{\mu}_{N-1}$, and the scenario where the IM is absent, \ie, when the true information states $\mu_{0},\cdots,\mu_{N-1}$ are used.
For this purpose, we denote by $\Tilde{J}_{N\pi\theta}$ the objective value of the DM when the IM is present and we focus on the quantity $|\Tilde{J}_{N\pi\theta}-J_{N\pi}|$.
Instead of elaborating the expression of this quantity based on the original POMDP defined in Section \ref{sec:model} using the random utility $U(\sum_{k=0}^{N-1}\beta^kc(X_k,Y_k,A_k))$ and the probability measures $\mathbb{P}^{\pi\theta}$ and $\mathbb{P}^{\pi}$, we resort to their representations under the equivalent MDP introduced in Section \ref{sec:solution:DP equations} using the information states.
That is, we consider
\begin{equation*}
\begin{aligned}
    |\Tilde{J}_{N\pi\theta}-J_{N\pi}|=& \Bigg|\int_{\mathcal{X}}\cdots\int_{\mathcal{X}}\int_{\mathcal{Y}}\int_{\mathbb{R}_+}U(s_N)\Tilde{\mu}_{N}(dy_N,ds_N)\lambda(dx_0)\cdots\lambda(dx_N)\\
    & \  \ - \int_{\mathcal{X}}\cdots\int_{\mathcal{X}}\int_{\mathcal{Y}}\int_{\mathbb{R}_+}U(s_N)\mu_{N}(dy_N,ds_N)\lambda(dx_0)\cdots\lambda(dx_N)\Bigg|.
\end{aligned}
\end{equation*}
Let $\Bar{c}:=||c(x,y,a)||_{\infty}$. The next result characterizes the impact of information manipulation.
}

\textcolor{black}{
\begin{theorem}
Let $\varepsilon_n:=\sup_{x\in X, y\in Y, a\in A}||P_n-Q(\cdot|x,y,a)||_{TV}$ denote the greatest possible information distortion resulting from manipulation $P_n$ for stage $n=1,\cdots,N-1$, 
and $\varepsilon_0:=||P_0-Q_0||_{TV}$ denote that for stage $0$.
Then, the performance deviation caused by information manipulation satisfies
\begin{equation}
    |\Tilde{J}_{N\pi\theta}-J_{N\pi}|\leq  \varepsilon_0 U(\Bar{c}\frac{1-\beta^N}{1-\beta})+\sum_{k=1}^{N-1} \left(\varepsilon_k U(\Bar{c}\frac{1-\beta^k}{1-\beta})\right).
    \label{eq:N stage performance deviation}
\end{equation}
\end{theorem}
\noindent{\bf Proof}.
We start with the following expression of the deviation between the manipulated and the true information states on a set $B\in\mathcal{B}(\mathcal{Y}\times\mathbb{R}_+)$ for $n=1,\cdots,N$:
\begin{equation*}
    \begin{aligned}
        |\Tilde{\mu}_n-\mu_n|(B)=&|\Tilde{\Psi}_n(x_{n-1},a_{n-1},x_n,\Tilde{\mu}_{n-1},\beta^{n-1})-\Psi_n(x_{n-1},a_{n-1},x_n,\mu_{n-1},\beta^{n-1})|(B) \\
        =& \int_{\mathcal{Y}}\int_{\mathbb{R}_+}\int_{B}\nu(dy_n)\delta_{s_{n-1}+\beta^{n-1}c(x_{n-1},y_{n-1},a_{n-1})}(ds_n)\\
        &\  \ \cdot \Big|p_n(x_n,y_n) \Tilde{\mu}_{n-1}(dy_{n-1},ds_{n-1})
        -
        q(x_n,y_n|x_{n-1},y_{n-1},a_{n-1})\mu_{n-1}(dy_{n-1},ds_{n-1})\Big|.
        \\
    \end{aligned}
\end{equation*}
Consider the auxiliary term $p_n(x_n,y_n)\mu_{n-1}(dy_{n-1},ds_{n-1})$ and apply triangle inequality, we arrive at $ |\Tilde{\mu}_n-\mu_n|(B)\leq e_n^1(B)+e_n^2(B)$, where the error terms $e_n^1$ and $e_n^2$ are defined as 
\begin{equation*}
    \begin{aligned}
        e_n^1(B):=&\int_{\mathcal{Y}}\int_{\mathbb{R}_+}\int_{B}\nu(dy_n)\delta_{s_{n-1}+\beta^{n-1}c(x_{n-1},y_{n-1},a_{n-1})}(ds_n)p_n(x_n,y_n)|\Tilde{\mu}_{n-1}-\mu_{n-1}|(dy_{n-1},ds_{n-1}),
        \\
        =&\int_{\mathbb{R}_+}\int_{B}p_n^X(x_n)\delta_{s_{n-1}+\beta^{n-1}c(x_{n-1},y_{n-1},a_{n-1})}(ds_n)|\Tilde{\mu}_{n-1}-\mu_{n-1}|(dy_{n-1},ds_{n-1}).
    \end{aligned}
\end{equation*}
and
\begin{equation*}
    \begin{aligned}
        e_n^2(B):=&\int_{\mathcal{Y}}\int_{\mathbb{R}_+}\int_{B}\nu(dy_n)\delta_{s_{n-1}+\beta^{n-1}c(x_{n-1},y_{n-1},a_{n-1})}(ds_n)\\
        &\  \ \cdot |p_n(x_n,y_n)-q(x_n,y_n|x_{n-1},y_{n-1},a_{n-1})|\mu_{n-1}(dy_{n-1},ds_{n-1}).
    \end{aligned}
\end{equation*}
The following inequality is a consequence of the above expression of $e_n^1$:
\begin{equation}
    \begin{aligned}
        &\int_{\mathcal{X}}\cdots\int_{\mathcal{X}}\int_{\mathcal{Y}}\int_{\mathbb{R}_+}U(s_n)e_n^1(dy_n,ds_n)\lambda(dx_0)\cdots\lambda(dx_n)\\
        =&\int_{\mathcal{X}}\cdots\int_{\mathcal{X}}\int_{\mathcal{Y}}\int_{\mathbb{R}_+} U(s_{n-1}+\beta^{n-1}c(x_{n-1},y_{n-1},a_{n-1}))|\Tilde{\mu}_{n-1}-\mu_{n-1}|(dy_{n-1},ds_{n-1})\lambda(x_0)\cdots\lambda(dx_{n-1})
        \\
        \leq&\int_{\mathcal{X}}\cdots\int_{\mathcal{X}}\int_{\mathcal{Y}}\int_{\mathbb{R}_+} U(s_{n-1}+\beta^{n-1}\Bar{c})|\Tilde{\mu}_{n-1}-\mu_{n-1}|(dy_{n-1},ds_{n-1})\lambda(x_0)\cdots\lambda(dx_{n-1}).
        \label{eq:proof:E[e_n^1]}
    \end{aligned}
\end{equation}
The error $e_n^2$ can be written as:
\begin{equation*}
    \begin{aligned}
        e_n^2(B)=&\int_{\mathcal{Y}}|p_n(x_n,y_n)-q(x_n,y_n|x_{n-1},y_{n-1},a_{n-1})|\nu(dy_n) \\
        & \  \ \cdot \int_{\mathbb{R}_+}\int_{B}\delta_{s_{n-1}+\beta^{n-1}c(x_{n-1},y_{n-1},a_{n-1})}(ds_n)\mu_{n-1}(dy_{n-1},ds_{n-1})\\
        =& \int_{\mathcal{Y}}|p_n(x_n,y_n)-q(x_n,y_n|x_{n-1},y_{n-1},a_{n-1})|\nu(dy_n) \\
        & \  \ \cdot \int_{B}\int_{\mathbb{R}_+}\delta_{s_{n-1}+\beta^{n-1}c(x_{n-1},y_{n-1},a_{n-1})}(ds_n) q(x_{n-1},y_{n-1}|x_{n-2},y_{n-2},a_{n-2})\nu(dy_{n-1}) \\
        & \  \ \cdot  \cdots \\
        & \  \ \cdot \int_{\mathbb{R}_+}\int_{\mathcal{Y}}\delta_{s_0+\beta^0 c(x_0,y_0,a_0)}(ds_0)\mu(dy_0,ds_0),
    \end{aligned}
\end{equation*}
where the last equality is obtained by applying (\ref{eq:info state operator}) repeatedly.
Then, we arrive at the following relation:
\begin{equation}
    \begin{aligned}
        &\int_{\mathcal{X}}\cdots\int_{\mathcal{X}}\int_{\mathcal{Y}}\int_{\mathbb{R}_+}U(s_n)e_n^2(dy_n,ds_n)\lambda(dx_0)\cdots\lambda(dx_n)\\
        \leq &      \int_{\mathcal{X}}\int_{\mathcal{Y}}\cdots\int_{\mathcal{X}}\int_{\mathcal{Y}}|p_n(x_n,y_n)-q(x_n,y_n|x_{n-1},y_{n-1},a_{n-1})|U\left(\sum_{k=0}^{n-1}\beta^kc(x_{k},y_{k},a_{k})\right)\\
        & \  \ \cdot 
        q(x_{n-1},y_{n-1}|x_{n-2},y_{n-2},a_{n-2})\cdots q(x_1,y_1|x_0,y_0,a_0)q_0^X(x_0)q_0^Y(y_0)\\
        &\  \ \cdot \nu(dy_{n})\cdots\nu(dy_0)\lambda(dx_{n})\cdots\lambda(dx_0)
        \\
        \leq&
        \varepsilon_n U(\sum_{k=0}^{n-1}\beta^k \Bar{c}).
    \end{aligned}
    \label{eq:proof:E[e_n^2]}
\end{equation}
Combining (\ref{eq:proof:E[e_n^1]}) and (\ref{eq:proof:E[e_n^2]}) and observing that $e_n^1$ and $e_n^2$ can be defined for $n=n-1,\cdots,1$, we obtain the following expression of the performance deviation:
\begin{equation*}
    \begin{aligned}
        |\Tilde{J}_{n\pi\theta}-J_{n\pi}|\leq&
        \int_{\mathcal{X}}\cdots\int_{\mathcal{X}}\int_{\mathcal{Y}}\int_{\mathbb{R}_+} U(s_{n-1}+\beta^{n-1}\Bar{c})|\Tilde{\mu}_{n-1}-\mu_{n-1}|(dy_{n-1},ds_{n-1})\lambda(x_0)\cdots\lambda(dx_{n-1})
        \\
        & \  \ + \varepsilon_n U(\sum_{k=0}^{n-1}\beta^k \Bar{c})
        \\
        \leq &\int_{\mathcal{X}}\cdots\int_{\mathcal{X}}\int_{\mathcal{Y}}\int_{\mathbb{R}_+} U(s_{n-2}+\Bar{c}\sum_{k=n-2}^{n-1}\beta^k)|\Tilde{\mu}_{n-2}-\mu_{n-2}|(dy_{n-2},ds_{n-2})\lambda(x_0)\cdots\lambda(dx_{n-2})\\
        & \  \ + \varepsilon_{n-1} U(\sum_{k=0}^{n-2}\beta^k \Bar{c})+ \varepsilon_n U(\sum_{k=0}^{n-1}\beta^k \Bar{c})
        \\
        \leq& \cdots \\
        \leq& U(\Bar{c}\sum_{k=0}^{n-1}\beta^k)||\Tilde{\mu}_0-\mu_0||_{TV}+\varepsilon_1 U(\Bar{c}\sum_{k=0}^{0}\beta^k)+\cdots+\varepsilon_{n-1} U( \Bar{c}\sum_{k=0}^{n-2}\beta^k)+ \varepsilon_n U(\Bar{c}\sum_{k=0}^{n-1}\beta^k )\\
        =& \varepsilon_0 U(\Bar{c}\frac{1-\beta^n}{1-\beta})+\sum_{k=1}^n \left(\varepsilon_k U(\Bar{c}\frac{1-\beta^k}{1-\beta})\right).
    \end{aligned}
\end{equation*}
Since the manipulations of the IM take place at stages $n=0,\cdots,N-1$, we arrive at the assertion in the theorem by observing that $e_N^2=0$.
\qed
}

\textcolor{black}{
The components on the right-hand side of the upper-bound (\ref{eq:N stage performance deviation}) have the following interpretations.
The first term of the bound corresponds to the impact on the performance caused by manipulated prior.
This source of impact lasts throughout the decision horizon, which can be observed from its dependence on the costs accumulated in all of the $N$ stages.
The second term of the upper-bound correspond to the impacts from the distorted beliefs.
Since a belief, or an information state, in our setting incorporates the (unnormalized) conditional distribution of the accumulated cost so far, the information manipulation at stage $k$ distorts the DM's conditional evaluation of the costs accumulated up to stage $k-1$.
}

The upper-bound (\ref{eq:N stage performance deviation}) is tight in the sense that the right-hand side of (\ref{eq:N stage performance deviation}) reduces to $0$ when the IM does not change the joint probability distributions, \ie, $\varepsilon_n=0$ for all $n=0,1,\cdots,N-1$.

The performance upper bound, (\ref{eq:N stage performance deviation}) can be used by the IM to estimate the minimum number of stages needed to drag the performance away from the original by a given amount.
To observe this, suppose that this target amount is $\mathfrak{g}>0$ in this adversarial scenario.
We assume for now that the effort of the IM in information distortion can achieve the amount $\Bar{\varepsilon}$ for all stages, $\ie$, $\varepsilon_n= \Bar{\varepsilon}$ for  $n=0,1,\cdots, N-1$.
Then, (\ref{eq:N stage performance deviation}) can be rewritten as $|\Tilde{J}_{N\pi\theta}-J_{N\pi}|\leq \Bar{\varepsilon} \sum_{j=1}^{N}U(\Bar{c}\frac{1-\beta^j}{1-\beta}) $.
Since the right-hand side of this relation is strictly increasing in the number of total stages $N$, we can solve for the smallest integer $N$ such that $\Bar{\varepsilon} \sum_{j=1}^{N}U(\Bar{c}\frac{1-\beta^j}{1-\beta}) \geq \mathfrak{g}$.
Then this integer has the interpretation of the minimal persistency requirement for the IM to achieve the manipulation goal $\mathfrak{g}$ if the effort $\Bar{\varepsilon}$ is feasible at each stage.

\textcolor{black}{
Staying stealthy is essential for the IM to successfully perform information manipulation.
One approach for the DM to verify the credibility of observed information is adopting statistical testing methods, which require a sufficient number of observation samples collected.
Hence, it is in general safer for the IM to complete the manipulation in a short period of time.
It is also feasible for the IM to constrain the distortion at each stage so that the detectability of the distributional shift is decreased.
The cost $\rho_n$ in the IM's objective function can be interpreted as a soft stealthy constraint, since the IM may not afford a design $P_n$ that significantly deviates from $Q$ given a large relative cost $\gamma$.
One may also consider using hard constraints on the IM's design problem at each stage, such as setting $\varepsilon_n\leq \Bar{\varepsilon}$ at each stage.
}

\textcolor{black}{
From (\ref{eq:N stage performance deviation}) we also observe that the risk preference of a DM, represented by $U$, influences the effectiveness of manipulation.
The IM may hope to adjust his manipulation strategies based on the knowledge of the DM's risk preference.
Consider the following thought experiment.
Suppose that the DM's prior knowledge of the system is not exploited by the IM, \ie, $\varepsilon_0=0$, and that the IM aims to distort the beliefs of the DM with his total information distortion effort limited by $\Bar{\varepsilon}$, \ie, $\varepsilon_1+\cdots+\varepsilon_{N-1}\leq \Bar{\varepsilon}$.
If one ignores the cost of performing manipulation, \ie, $\rho_k=0$ for $k=1,\cdots,N-1$, then the optimal strategy for the IM is to allocate $\Bar{\varepsilon}$ only to stage $N-1$.
When the cost of performing manipulation is positive, there is a need for the IM to tradeoff between the impact and the cost of allocating the distortion effort to a single stage.
If the DM is risk-seeking, \ie, $U$ is concave, the marginal increase in the DM's perceived cost $U(\sum_{n=0}^{N-1}\beta^k c(x_k,y_k,a_k))$ decreases as the instantaneous cost $c$ accumulates.
Then, it is more likely that the impacts of allocating the same amount of distortion effort to nearby stages are similar.
Consequently, the IM would prefer to distribute his manipulation effort to avoid large costs of performing manipulation.
However, if the DM is risk-averse, \ie, $U$ is convex, the allocation result would be more concentrated on later stages.
The reason lies in that the marginal impact increase is more likely to be greater than the marginal manipulation cost increase due to convexity.
}

\section{Concluding Remarks}
\label{sec:conclusion}
Motivated by the informational disadvantages of decision-makers in partially observable systems, we have proposed DIMG to study the generation and the impact of misinformation from the perspective of an IM.
We have built our framework on a generic POMDP to understand how DMs interpret misleading information and how these interpretations influence decisions in a dynamic setting.

The solution to the DIMG relies on the information state approach to solving POMDPs. 
We have first presented the technique of using unnormalized joint conditional densities of the hidden state and the accumulated cost to represent a risk-sensitive DM's objective function.
This way of reformulating the POMDP has simplified its dynamic programming equations.
By identifying the DIMG at each stage and leveraging the dynamic programming equations, we have shown the existence of the equilibria of information distortion under both the ex ante and the interim information patterns at each stage. 
The stage-wise consistency constraint for restricting the IM's power connects the stage problems and leads to the solution to the IM's manipulation problem.
Examples have been investigated to show the reformulation of a DIMG into a linear programming problem and to elaborate on the relation between the maneuverability of information and the properties of the underlying system.
We have also derived a tight upper bound on the impact of information manipulation on the performances of the system.
Our stability results have shown that, seen from the final stage, IM's influence on the system is bounded by a quantity which consists of instantaneous informational differences at certain stages together with the expected stage-wise amplifications at their subsequent stages.
DMs with high-risk aversions are more likely to resist the impacts of information manipulation.

We conclude the paper with the following applications.
\paragraph{Defensive gaslighting}
Recent years have witnessed a significant increase in computer network failures and cyber system interruptions caused by cyber risks.
Proactive defense mechanisms constitute a class of system protection techniques aiming at enhancing cyber system robustness actively before the deployments of cyber attacks.
The DIMG introduced in this paper is applicable to defensive gaslighting (e.g., \cite{liu2023impact}), which belongs to the class of cyber deception methods (see \cite{wang2018cyber}) and targets at informationally deceiving adversaries and attackers so that they either cease attack attempts or miss the essential nodes in a computer network system.

The attacker's planning of adversarial efforts can be modeled by a POMDP (e.g., \cite{anwar2022honeypot}).
Using the example containing discrete state spaces shown in Section \ref{sec:example} is one possibility, which indicates that , the unobservable state of the POMDP describes whether a subset of nodes (target of the attacker) in a computer network is compromisable, \ie, $\mathcal{Y}=\{``C", ``U"\}$ with $``C"$ meaning ``Compromisable" and $``U"$ meaning ``Uncompromisable".
The compromisability of the target is reflected partially by its characteristics during operation, such as input/output queues and communication links with other nodes.
One can model these characteristics in an abstract perspective and encapsulate them as four observable states $\mathcal{X}=\{``WC", ``SC", ``WU", ``SU"\}$, where $``W"$ means ``Weakly" and $``S"$ means ``Strongly".
Depending on system structures and defensive investments of nodes, the transition kernel may vary.
The setting in Section \ref{sec:example} is one possibility, which can be interpreted as that ``Compromisable" nodes can be observed as ``Strongly Compromisable" and ``Weakly Compromisable", and that ``Uncompromisable" nodes can be observed as ``Strongly Uncompromisable" and ``Weakly Uncompromisable".
The network defender plays the role of the IM in the DIMG.
The defender aims to manipulate the characteristics of the nodes so that key functional nodes, such as data storage centers and nodes with massive connections or high communication frequency, are excluded from the set of direct attack targets of the attacker.
The defender can achieve this goal by resorting to the linear programming problem discussed in Section \ref{sec:example}, which can be solved within a reasonable amount of time using existing optimization software.
Ideal outcomes of the information manipulation include making honeypot nodes in the network to be observed as ``Strongly Compromisable" so that the attacker's effort is wasted and presenting vulnerable key nodes as ``Strongly Uncompromisable" so that they function properly.

\paragraph{Objective-driven ambiguity set construction.}
Ambiguity sets are crucial for optimal decision making in distributionally robust optimization problems.
Popular approaches for constructing ambiguity sets include, for example, the Wasserstein ball approach (e.g., \cite{esfahani2015data}) and the moment-type approach (e.g., \cite{guo2022robust}).
Nominal distributions of uncertainty are essential in the first approach.
The nominal distribution can be predetermined by the DM using prior knowledge about the decision environment or represented by the empirical distribution derived from available data.
The solutions to the IM's problem in our framework are also a meaningful candidate for the nominal distribution.
If they are the centers of the ambiguity sets, then, these sets can be considered as driven by the IM's objective, since the origins of the ambiguities arise from the presence of the IM.

\paragraph{Information manipulation in principal-agent problems}
Principal-agent problems have been studied extensively due to their broad application domains.
The relation between the IM and the DM in the DIMG resembles that of the relation between the principal and the agent.
Therefore, one can apply information manipulation to provide the principal with an additional maneuver to enhance her design.

One application domain of principal-agent problems is congestion control.
Efficient transportation strategies for congestion mitigation are fundamental in building smart cities.
With the rapid advancements of informational technologies, smart routing and congestion mitigation methods become more available to drivers and more adaptable to various traffic circumstances.
Drivers on the road observe little traffic information.
They need to resort to routing recommendations, which serve as partial observations about the global traffic information that influences their travels from origins to destinations.
Authorities can leverage the DIMG to dynamically generate recommendations aiming at encouraging users to adopt routing strategies that are socially optimal (see \cite{wardrop1952road}).

Another principal-agent relation where one can involve information manipulation is contract design.
In a class of contract problems, the distribution of the state variable which influences the principal's payoff depends on the agent's action.
This relation is naturally analogous to the fact that the transition law of the unobservable state in our framework depends on the DM's actions.
When a dynamic contract is considered, the principal's design drives the agent's actions so that the distribution of the unobservable states favors the principal.
At the same time, the IM can also play with the distributions of the observable states, based on which the DM takes actions.
Information maneuvering serves as an alternative to designing contracts with incentive compatibility (e.g., \cite{liu2024stackelberg}).
This effect can potentially reduce the gap between the first-best and second-best contracts, resulting in an increased practicality of the first-best contract. 
Disregarding the concerns of legality and fairness, information manipulation can further increase the principal's revenue and be computationally tractable in practice.
The risk preference influence problem considered in \cite{liu2024stackelberg} can also improve the effect of information manipulation. 
The reason lies in that, as we have mentioned at the end of Section \ref{sec:impacts}, the impact of information manipulation is dependent on the degree of risk aversion of a DM.

\bibliographystyle{unsrt}  
\bibliography{references}  

\nocite{*}

\end{document}